\documentclass[a4paper,oneside,11pt]{article}

\usepackage{amsfonts,color}

\newtheorem{dfn}{Definition}[section]
\newtheorem{tw}[dfn]{Theorem}
\newtheorem{prop}[dfn]{Proposition}
\newtheorem{rem}[dfn]{Remark}
\newtheorem{ex}[dfn]{Example}
\newtheorem{lem}[dfn]{Lemma}

\newtheorem{cor}[dfn]{Corollary}
\usepackage{amssymb} 
\usepackage{amsmath}
\numberwithin{equation}{section}

 \global\long\def\sbr#1{\left[ #1\right] }
 \global\long\def\cbr#1{\left\{  #1\right\}  }
 
 \global\long\def\rbr#1{\left(#1\right)}

 \global\long\def\R{\mathbb{R}}
\global\long\def\S{\mathbb{S}}

 \global\long\def\dd#1{\textnormal{d}#1}

 \global\long\def\ra{\rightarrow}
 
 \global\long\def\ns{\infty}

\usepackage{anysize}
\usepackage{array}
\usepackage{enumerate}

\title{\bf CIR equations with multivariate L\'evy noise}

\author{Micha\l \ Barski, Rafa\l \ \L ochowski}

\begin{document}

\maketitle

\begin{abstract}
The paper is devoted to the study of the short rate equation of the form
\begin{gather*}
\dd R(t)=F(R(t))\dd t+\sum_{i=1}^{d}G_i(R(t-))\dd Z_i(t), \quad R(0)=x\geq 0,\quad t>0,
\end{gather*}
with deterministic functions $F,G_1,...,G_d$ and a multivariate L\'evy process $Z=(Z_1,...,Z_d)$. The equation is supposed 
to have a nonnegative solution which generates an affine term structure model. Two classes of noise are considered. In the first one the coordinates of Z are independent processes with regularly varying Laplace exponents. In the second class Z is a spherical processes, which means that its L\'evy measure has a similar structure as that of a stable process, but with radial part of a general form. For both classes a precise form of the short rate generator is characterized. Under mild assumptions it is shown that any equation of the considered type has the same solution as the equation driven by a L\'evy process with independent stable coordinates.

The paper generalizes the classical results on the Cox-Ingersoll-Ross (CIR) model, \cite{CIR}, as well as on its extended version from \cite{BarskiZabczykCIR} and \cite{BarskiZabczyk} where $Z$ is a one-dimensional L\'evy process.
\end{abstract}

\tableofcontents

\section{Introduction}
This paper is concerned with a stochastic equation of the form 
\begin{gather}\label{rownanie 1}
\dd R(t)=F(R(t))\dd t+\sum_{i=1}^{d}G_i(R(t-))\dd Z_i(t), \quad R(0)=x,\quad t>0,
\end{gather}
where $F$, $\{G_i\}_{i=1,2,...,d}$ are deterministic functions, $Z_i(t), i=1,2,...,d$, are L\'evy processes and martingales,
 $x$ is a nonnegative constant. A solution $R(t), t\geq 0$, if nonnegative, will be identified here with the short rate process, so
it defines the bank account process by
$$
B(t):=e^{\int_{0}^{t}R(s)ds}, \quad t\geq 0.
$$
Related to the savings account are zero coupon bonds. Their prices form a family of stochastic processes
$P(t,T), t\in [0,T]$, parametrized by their maturity times $T\geq 0$.  The price of a bond with maturity $T$ at time $T$ is equal to its
nominal value, typically assumed, also here, to be $1$, that is $P(T,T)=1$. The family of bond prices is supposed to have the {\it affine structure}, which means that 
\begin{gather}\label{affine model}
P(t,T)=e^{-A(T-t)-B(T-t) R(t)}, \quad 0\leq t\leq T,
\end{gather}
for some smooth deterministic functions $A$, $B$. Hence, the only source of randomness in the affine model \eqref{affine model} is the short rate process $R$ given by \eqref{rownanie 1}.  As the resulting market constituted by $(B(t), \{P(t,T)\}_{T\geq 0})$
should exclude arbitrage, the discounted bond prices
$$
\hat{P}(t,T):=B^{-1}(t)P(t,T)=e^{-\int_{0}^{t}R(s)ds-A(T-t)-B(T-t)R(t)}, \quad 0\leq t\leq T,
$$
are supposed to be local martingales for each $T\geq 0$. This requirement affects in fact our starting equation. Thus  
the functions $F$, $\{G_i\}_{i=1,...,d}$ and the noise $Z=(Z_1,...,Z_d)$ should be chosen such that $\eqref{rownanie 1}$ has a nonnegative solution with any $x\geq 0$ and such that, for some functions $A,B$ and each $T\geq 0$, \ $\hat{P}(t,T)$ is a local martingale on $[0,T]$. If this is the case, \eqref{rownanie 1} will be called to {\it generate an affine model} or to be a {\it generating equation}, for short.

In the case  when $Z=W$ is a real-valued Wiener process, the only generating equation is the classical CIR equation
\begin{gather}\label{CIR equation}
\dd R(t)=(aR(t)+b)\dd t+C\sqrt{R(t)}\dd W(t), 
\end{gather}
with $a\in\mathbb{R}$, $b,C\geq 0$, due to Cox, Ingersoll, Ross, see \cite{CIR}. The 
case with a general one-dimensional L\'evy process $Z$ was studied in \cite{BarskiZabczykCIR}, \cite{BarskiZabczyk} and
\cite{BarskiZabczykArxiv} with the following conclusion. If the variation of $Z$ is infinite and $G \not\equiv 0$, then $Z$ must be an $\alpha$-stable process with index $\alpha\in(1,2]$, with either positive or negative jumps only,  and \eqref{rownanie 1} has the form
\begin{gather}\label{CIR equation geenralized}
\dd R(t)=(aR(t)+b)\dd t+C\cdot R(t)^{{1}/{\alpha}}\dd Z(t),
\end{gather}
with $a\in\mathbb{R}, b\geq 0$ and $C$ such that it has the same sign as the jumps of $Z$. Clearly, for $\alpha=2$ equation \eqref{CIR equation geenralized} becomes \eqref{CIR equation}. If $Z$ is of finite variation then the noise enters \eqref{rownanie 1} in the additive way, that is 
\begin{gather}\label{Vasicek equation geenralized}
\dd R(t)=(aR(t)+b)\dd t+C \ \dd Z(t).
\end{gather}
Here $Z$ can be chosen as an arbitrary process with positive jumps, $a\in\mathbb{R}, C\geq 0$ and 
$$
b\geq C \int_{0}^{+\infty}y \ \nu(\dd y),
$$
where $\nu(\dd y)$ stands for the L\'evy measure of $Z$.  The variation of $Z$ is finite, so is the right side above.
Recall, \eqref{Vasicek equation geenralized} with $Z$ being a Wiener process is the well known Vasi\v cek equation, see \cite{Vasicek}. Then the short rate is a Gaussian process, hence it takes negative values with positive probability. 
This drawback is eliminated by the jump version of  the Vasi\v cek equation \eqref{Vasicek equation geenralized}.

This paper is devoted to the equation \eqref{rownanie 1} with $d>1$. The multidimensional setting makes the study of equation 
\eqref{rownanie 1} more complicated. The reason is that, unlikely as in the case $d=1$, different generating equations may have identical solutions in the sense that the solutions' generators are the same. Our first goal is to characterize the class of generators of solutions of generating equations and the second goal is to construct, for each element of this class,  a related specific equation. In this way any short rate process given by \eqref{rownanie 1} which generates an affine model becomes representable by a tractable equation. This approach seems to be useful for future applications. 

Our solution of the problem is based on, rather abstract, result of Filipovi\'c \cite{FilipovicATS}, characterizing generators of a general Markovian non-negative short rate process. The contribution of this paper is making this characterization concrete for two classes of L\'evy processes.  In the first class the coordinates of the noise 
$$
Z_1(t),Z_2(t),...,Z_d(t), \quad t\geq 0,
$$ 
are independent L\'evy processes being martingales of infinite variation. Their Laplace exponents are assumed to vary regularly at zero. We show that the solution of any generating equation with such a noise is the same as the solution of the equation
$$
\dd R(t)=(a R(t)+b)\dd t+\sum_{k=1}^{g} d_k R(t-)^{1/\alpha_k} \dd Z^{\alpha_k}_k(t),
$$
where $1\leq g\leq d$, $a\in\mathbb{R}$, $b\geq 0$, $d_k>0$, $2\geq \alpha_1>...>\alpha_g>1$, and $Z_k^{\alpha_k}$ is a stable process with index $\alpha_k$. The second class consists of 
{\it spherical L\'evy processes}. We call a process $Z(t):=(Z_1(t),...,Z_2(t))$ spherical if its 
 L\'evy measure $\nu(\dd y)$ admits the following representation
\begin{equation} \label{spherical}
\nu(A)=\int_{\S^{d-1}}\lambda(\dd \xi)\int_{0}^{+\infty}\mathbf{1}_{A}(r\xi)\gamma(\dd r), \quad A\in\mathcal{B}(\mathbb{R}^d).
\end{equation}
 Here $\S^{d-1}$ is a unit sphere in $\mathbb{R}^d$, $\lambda(\dd \xi)$ is a finite measure on $\S^{d-1}$\ called a {\it spherical part} of $\nu$, $\gamma(\dd r)$ is a L\'evy measure on $(0,+\infty)$ called a {\it radial part} of $\nu$. One can see that 
on each half-line in $\mathbb{R}^d$ starting from the origin, the L\'evy measure is given in the same way by the radial part,
up to multiplication by a nonnegative constant.  An important example of a radial measure satisfying 
\eqref{warunek na miare radialna a}-\eqref{warunek na miare radialna b} is
\begin{gather}\label{radialna stabilna}
\gamma(\dd r)=\frac{1}{r^{1+\alpha}}\dd r, \quad \alpha\in(1,2).
\end{gather}
Given this measure and any finite measure $\lambda(d\xi)$,  the formula  
\eqref{spherical} corresponds to a stable process with index $\alpha\in(1,2)$. This process has no Wiener part. The $2$-stable process is the Wiener process. We prove, under mild conditions, that the solution of any generating equation with spherical noise is the same as the solution of equation \eqref{CIR equation geenralized}.

Our results for each of the classes introduced above generalize the one dimensional results from \cite{BarskiZabczykCIR}, \cite{BarskiZabczyk} and \cite{BarskiZabczykArxiv}.

The structure of the paper is as follows. In Section \ref{section Formulation of the problem} we introduce the probabilistic setting
for the equation \eqref{rownanie 1} and present a properly adapted version of the result from  \cite{FilipovicATS} characterizing the generator of a generating equation. In particular, we point out here the role of the projections of $Z$ along $G$, meant as processes $\sum_{i=1}^{d}G_i(x)Z_i(t)$, $x\geq 0$, for the generator of $R$. Using examples highlighting the differences between the one- and multidimensional case we justify the form of problem-stating described above. Section \ref{section Noise with independent coordinates} is concerned with equation \eqref{rownanie 1} driven by $Z$ with independent coordinates. The main result here is Theorem \ref{TwNiez}. Regularly varying Laplace exponents are described in terms of the L\'evy measure in Subsection \ref{sec slowly varying Laplace exp.}. We also describe all generating equations in the case $d=2$ in Subsection \ref{section Generalized CIR equations on a plane} and provide an example showing the non-uniqueness of a generating equation when $d=3$ in Subsection \ref{section Example in higher dimensions}. The case when $Z$ is spherical is presented in Section \ref{section Spherical Levy noise}.
The main result here is Theorem \ref{tw_sferyczne}. Its proof, presented in Subsection \ref{proof_of_tw_sf}, requires a sequence of auxiliary results contained in Subsection \ref{Auxilliary results}.

\section{Preliminaries}\label{section Formulation of the problem}

The problem of description of generating equations \eqref{rownanie 1} in the multidimensional case will be handled in a different way than in the one-dimensional case. Basing on Proposition \ref{prop wstepny} in Subsection \ref{section Projections of the noise} and examples in Subsection \ref{sec Examples} we explain here the formulation of the problem studied in the sequel.

\subsection{Setup for the equation}
Using  the scalar product $\langle\cdot,\cdot\rangle$ in $\mathbb{R}^d$ we write \eqref{rownanie 1} in the short form
\begin{gather}\label{rownanie 2}
\dd R(t)=F(R(t))\dd t+\langle G(R(t-)),\dd Z(t)\rangle, \quad R(0)=x\geq 0, \qquad t>0,
\end{gather}
where $F:[0,+\infty)\longrightarrow \mathbb{R}$, $G:=(G_1,G_2,...,G_d):[0,+\infty)\longrightarrow\mathbb{R}^d$ and 
 $Z:=(Z_1,Z_2,...,Z_d)$  is a L\'evy process in $\mathbb{R}^d$ with the characteristic triplet $(a,Q,\nu(\dd y))$. Recall, $a\in\mathbb{R}^d$ describes the drift part of $Z$, $Q$ is a non-negative, symmetric, $d\times d$ covariance matrix, characterizing the coordinates' covariance of the Wiener part $W$ of $Z$, and $\nu(\dd y)$ is a measure on $\mathbb{R}^d\setminus\{0\}$ satisfying
\begin{gather}\label{warunek na miare Levyego 1}
\int_{\mathbb{R}^d}(\mid y\mid^2\wedge \ 1)\ \nu(\dd y)<+\infty,
\end{gather}
describing the jumps of $Z$ and called the L\'evy measure of $Z$. Recall, $Z$ admits a representation as a sum of four independent processes of the form
\begin{gather}\label{LevyIto}
Z(t)=at +W(t)+\int_{0}^{t}\int_{\{\mid y\mid\leq 1\}}y\tilde{\pi}(\dd s,\dd y)+\int_{0}^{t}\int_{\{\mid y\mid> 1\}}y\pi(\dd s,\dd y),
\end{gather}
called the L\'evy-It\^o decomposition of $Z$. Above $\pi(\dd s,\dd y)$ and $\tilde{\pi}(\dd s,\dd y):=\pi(\dd s,\dd y)-\dd s \nu(\dd y)$ stand for the jump measure and the compensated jump measure of $Z$, respectively.
We consider the case when $Z$ is a martingale and call it a L\'evy martingale for short. Its drift and the L\'evy measure are such that 
\begin{gather}\label{warunek na miare Levyego 2}
\int_{\mid y\mid>1}\mid y\mid\ \nu(\dd y)<+\infty, \quad a+\int_{\mid y\mid>1}y \ \nu(\dd y)=0.
\end{gather}
Consequently, the characteristic triplet of $Z$ is 
\begin{gather}\label{chrakterystyki Z}
\left(-\int_{\mid y\mid>1}y \ \nu(\dd y), \ Q, \ \nu(\dd y)\right),
\end{gather}
and \eqref{LevyIto} takes the form
$$
Z(t)=W(t)+X(t), \qquad X(t):=\int_{0}^{t}\int_{\mathbb{R}^d}y \ \tilde{\pi}(\dd s,\dd y), \quad t\geq 0,
$$
where $W$ and $X$ are independent. The martingale $X$ will be called the jump part of $Z$. Its Laplace exponent $J_{\nu}$,  defined by the equation
\begin{equation} 
\mathbb{E}\sbr{e^{-\langle \lambda,X(t)\rangle}}=e^{tJ_{\nu}(\lambda)}, 
\end{equation}
has the  following representation
\begin{equation} \label{Jdef}
J_\nu(\lambda)=\int_{\mathbb{R}^d}(e^{-\langle\lambda,y\rangle}-1+\langle\lambda,y\rangle)\nu(\dd y),
\end{equation}
and is finite for $\lambda\in\mathbb{R}^d$ satisfying
$$
\int_{\mid y\mid>1}e^{-\langle \lambda,y\rangle}\nu(\dd y)<+\infty.
$$
By the independence of $X$ and $W$ we see that
$$
\mathbb{E}\sbr{e^{-\langle \lambda,Z(t)\rangle}}=\mathbb{E}\sbr{e^{-\langle \lambda,W(t)\rangle}}\cdot\mathbb{E}\sbr{e^{-\langle \lambda,X(t)\rangle}},
$$
so the Laplace exponent $J_Z$ of $Z$ equals
\begin{equation} \label{LaplaceZ}
J_Z(\lambda)={{\frac{1}{2}\langle Q\lambda,\lambda\rangle+J_\nu(\lambda)}}.
\end{equation}

\subsection{Projections of the noise and problem formulation}\label{section Projections of the noise}
For the function $G$ and the process $Z$ we consider the {\it projections} of $Z$ along $G$ given by
\begin{gather}\label{projection of Z}
Z^{G(x)}(t):=\langle G(x), Z(t)\rangle, \quad t\geq 0.
\end{gather}
For any $x\geq 0$, $Z^{G(x)}$ is a real-valued L\'evy martingale. It follows from the identity
$$
\mathbb{E}\sbr{e^{- \gamma\cdot Z^{G(x)}(t)}}=\mathbb{E}\sbr{e^{-\langle \gamma G(x), Z(t)\rangle}}, \quad \gamma\in\mathbb{R}, 
$$
and \eqref{LaplaceZ} that the Laplace exponent of $Z^{G(x)}$ equals 
\begin{gather}\label{Laplace ZG do uproszczenia}
J_{Z^{G(x)}}(\gamma)=J_Z(\gamma G(x))=\frac{1}{2}\gamma^2\langle Q G(x),G(x)\rangle+\int_{\mid y\mid>0}\left(e^{-\gamma \langle G(x),y\rangle}-1+\gamma\langle G(x),y\rangle\right)\nu(\dd y).
\end{gather}
Formula \eqref{Laplace ZG do uproszczenia} can be written in a simpler form by using the L\'evy measure $\nu_{G(x)}(\dd v)$ of $Z^{G(x)}$, which is the {\it image} 
of the L\'evy measure $\nu(dy)$ under the linear transformation $y\mapsto \langle G(x), y\rangle$. This measure will be denoted by
$\nu_{G(x)}(\dd v)$ and is given by 
$$
\nu_{G(x)}(A):=\nu \{y \in \R^d: \langle G(x),y\rangle\in A \} , \quad A\in\mathcal{B}(\mathbb{R}).
$$
Then we obtain from \eqref{Laplace ZG do uproszczenia} that
\begin{gather}\label{Laplace ZG}
J_{Z^{G(x)}}(\gamma)=\frac{1}{2}\gamma^2\langle Q G(x),G(x)\rangle+\int_{\mid v\mid>0}\left(e^{-\gamma v}-1+\gamma v \right)\nu_{G(x)}(\dd v).
\end{gather}
Thus the characteristic triplet of the projection $Z^{G(x)}$ has the form
\begin{gather}\label{charakterystyki rzutu}
\left(-\int_{\mid v\mid>1}y \ \nu_{G(x)}(\dd v), \ \langle Q G(x),G(x)\rangle, \ \nu_{G(x)}(\dd v)\mid_{v\neq 0}\right).
\end{gather}
Above we used the restriction $\nu_{G(x)}(\dd v)\mid_{v\neq 0}$ by cutting off zero which may be an atom of $\nu_{G(x)}(\dd v)$.

In Proposition \ref{prop wstepny} below we provide a preliminary characterization of equations \eqref{rownanie 2} generating affine models. The central role here is played by the law of $Z^{G}$. The result is deduced from Theorem 5.3 in \cite{FilipovicATS}, where the generator of a general non-negative Markovian  short rate process for affine models was characterized. The result settles a starting point for proving the main results of the paper.
\begin{prop}\label{prop wstepny} 
Let $Z$ be a L\'evy martingale with characteristic triplet  \eqref{chrakterystyki Z}, $Z^{G(x)}$ be its projection \eqref{projection of Z} and $\nu_{G(x)}(\dd v)$ be the L\'evy measure of $Z^{G(x)}$.
\begin{enumerate}[(A)] 
\item Then equation \eqref{rownanie 1} generates an affine model if and only if the following conditions are satisfied
\begin{enumerate}[a)]
\item For each $x \ge 0$ the support of $\nu_{G(x)}$ is contained in $[0, +\ns)$ which means that $Z^{G(x)}$ has positive jumps only, i.e.  for each $t\geq 0$, with probability one,
\begin{gather}\label{Z^G positive jumps}
\triangle Z^{G(x)}(t):=Z^{G(x)}(t)-Z^{G(x)}(t-)=\langle G(x), \triangle Z(t)\rangle\geq 0.
\end{gather}
\item The jump part of $Z^{G(0)}$ has finite variation, i.e.
\begin{gather}\label{nu G0 finite variation}
\int_{(0,+\infty)}v \ \nu_{G(0)}(\dd v)<+\infty.
\end{gather}
\item The characteristic triplet \eqref{charakterystyki rzutu} of $Z^{G(x)}$ is linear in $x$, i.e.
\begin{align}\label{mult. CIR condition}
\frac{1}{2}\langle Q G(x), G(x)\rangle&=cx, \quad x\geq 0,\\[1ex]\label{rozklad nu G(x)}
\nu_{G(x)}(\dd v)\mid_{(0,+\infty)}&=\nu_{G(0)}(\dd v)\mid_{(0,+\infty)}+x\mu(\dd v), \quad x\geq 0,
\end{align}
for some $c\geq 0$ and a measure  $\mu(\dd v) \ \text{on} \ (0,+\infty) \ \text{satisfying}$
\begin{gather}\label{war calkowe na mu}
\int_{(0,+\infty)}(v \wedge v^2)\mu(\dd v)<+\infty.
\end{gather}
\item The function $F$ is affine, i.e.
\begin{gather}\label{linear drift}
F(x)=ax+b, \ \text{where} \ a\in\mathbb{R}, \ b\geq\int_{(1,+\infty)}(v-1)\nu_{G(0)}(\dd v) .
\end{gather}
\end{enumerate}
\item Equation \eqref{rownanie 1} generates an affine model if and only if the generator of $R$ is given by
\begin{align}\label{generator R w tw}\nonumber
\mathcal{A}f(x)=cx f^{\prime\prime}(x)&+\Big[ax +b+\int_{(1,+\infty)}(1 -v)\{\nu_{G(0)}(\dd v)+x\mu(\dd v)\}\Big]f^{\prime}(x)\\[1ex]
&+\int_{(0,+\infty)}[f(x+v)-f(x)-f^{\prime}(x)(1\wedge v)]\{\nu_{G(0)}(\dd v)+x\mu(\dd v)\}.
\end{align}
for $f\in\mathcal{L}(\Lambda)\cup C_c^2(\mathbb{R}_{+})$, where 
$\mathcal{L}(\Lambda)$ is the linear hull of $\Lambda:=\{f_\lambda:=e^{-\lambda x}, \lambda\in(0,+\infty)\}$
and $C_c^2(\mathbb{R}_{+})$ stands for the set of twice continuously differentiable functions with compact support in $[0,+\infty)$. The constants $a,b,c$ and the measures $\nu_{G(0)}(\dd v), \mu(\dd v)$ are those from the part (A).
\end{enumerate}
\end{prop}

The poof of Proposition \ref{prop wstepny} is postponed to Appendix. 

In the sequel we use an equivalent formulation of  \eqref{mult. CIR condition}-\eqref{rozklad nu G(x)} with the use of Laplace exponents. Taking into account \eqref{Laplace ZG} we obtain the following.

\begin{rem}\label{rem warunki w eksp. Laplacea}
	The conditions \eqref{mult. CIR condition} and \eqref{rozklad nu G(x)}  are equivalent to
	\begin{gather}\label{war na exp Laplaca}
	J_{Z^{G(x)}}(b)=J_Z(bG(x))=cb^2x+J_{\nu_{G(0)}}(b)+x J_{\mu}(b), \quad b,x\geq 0,
	\end{gather}
	where
	\begin{gather}\label{def Jmu i Jnu0}
	J_{\mu}(b):=\int_{0}^{+\infty}(e^{-bv}-1+bv)\mu(\dd v), \quad J_{\nu_{G(0)}}(b):=\int_{0}^{+\infty}(e^{-bv}-1+bv)\nu_{G(0)}(\dd v).
	\end{gather}
\end{rem}

The part $(A)$ states that \eqref{rownanie 1} generates an affine model if $F$ is affine and 
the projections $Z^{G(x)}, x\geq 0$, have characteristic triplets characterized by a constant 
$c\geq 0$ carrying information on the presence of the Wiener part and two measures 
$\nu_{G(0)}(\dd v)$, $\mu(\dd v)$ describing jumps. In view of part $(B)$, the triplet $(c,\nu_{G(0)}(\dd v),\mu(\dd v))$
satisfying \eqref{nu G0 finite variation}-\eqref{war calkowe na mu}, together with $F$, determine the generator of the short rate process. A pair $(G,Z)$ for which the projections $Z^{G(x)}$ satisfy \eqref{nu G0 finite variation}-\eqref{war calkowe na mu}
will be called {\it a generating pair}. As $F$ is of a simple affine form, the essential issue is to characterize the measures $\nu_{G(0)}(\dd v)$ and $\mu(\dd v)$. In the one dimensional case the measures turn out to be such that either
$$
\bullet \quad  \nu_{G(0)}(\dd v)- \text{ is any measure on $(0,+\infty)$ of finite variation} \ \text{and} \ \mu(\dd v)\equiv 0,
$$
or
$$
\bullet \quad \nu_{G(0)}(\dd v)\equiv 0 \ \text{and} \ \mu(\dd v)=\frac{1}{v^{1+\alpha}}dv, v\geq 0,  \ \alpha\in(1,2), \ \text{i.e. $\mu(\dd v)$ is $\alpha$-stable}.
$$ 
Moreover, for each of the two situations above, there exists a unique, up to multiplicative constants, corresponding generating pair $(G,Z)$. For instance,
if $\mu(\dd v)$ is $\alpha$-stable then $Z$ is also $\alpha$-stable and $G(x)=C x^{{1}/{\alpha}}$, with some $C>0$.
This means that the triplet $(c,\nu_{G(0)}(\dd v),\mu(\dd v))$ corresponds to a unique equation \eqref{rownanie 2}, up to the choice of $F$.  The one-dimensional equations mentioned in Introduction can be characterized in terms of 
$(c,\nu_{G(0)}(\dd v),\mu(\dd v))$  in the following way.
\begin{enumerate} [a)]
\item $c>0, \ \nu_{G(0)}(\dd v)\equiv 0, \ \mu(\dd v)\equiv 0$; \\[1ex]
This case corresponds to the classical CIR equation \eqref{CIR equation}.
\item $c=0, \ \nu_{G(0)}(\dd v)\equiv0, \ \mu(\dd v)- \text{$\alpha$-stable}, \ \alpha\in(1,2)$;\\[1ex]
In this case \eqref{rownanie 2} becomes the generalized CIR equation with $\alpha$-stable noise \eqref{CIR equation geenralized}.
\item $c=0, \ \nu_{G(0)}(\dd v)- \text{any measure on $(0,+\infty)$ of finite variation}, \ \mu(\dd v)\equiv 0$;\\[1ex]
Here \eqref{rownanie 2} becomes the generalized Vasi\v cek equation \eqref{Vasicek equation geenralized}.
\end{enumerate}

In the case $d>1$ one should not expect a one to one correspondence between the triplets $(c,\nu_{G(0)}(\dd v),\mu(\dd v))$ and the generating equations \eqref{rownanie 2}. The reason is that the distribution of the product 
$\langle G(x), Z(t)\rangle$ does not determine the pair $(G,Z)$ in a unique way. If $(G,Z)$ and $(G^\prime,Z^\prime)$ are generating pairs satisfying \eqref{nu G0 finite variation}-\eqref{war calkowe na mu} with an identical triplet $(c,\nu_{G(0)}(\dd v),\mu(\dd v))$, then it follows from part (B) that the corresponding equations  \eqref{rownanie 2} have solutions with the same generator. Our illustrating examples in Section \ref{sec Examples} show a couple of different equations all providing the same short rate $R$. Furthermore, it turns out that even for a fixed process $Z$, the function $G$ in the generating pair $(G,Z)$ does not need to be unique. For this reason we focus in this paper on the characterization of possible laws of projections $Z^G$. 
As in the classes of L\'evy processes under our consideration $G(0)=0$, hence $\nu_{G(0)}(\dd v)$ vanishes, the goal is to determine the measure $\mu(\dd v)$ in the multidimensional case. Next, having such laws we provide corresponding generating equations which are of tractable form.  

\vskip2ex

\subsection{Examples}\label{sec Examples}

We present a couple of examples of generating pairs $(G,Z)$ such that the related projections $Z^{G(x)}$ satisfy 
conditions \eqref{nu G0 finite variation}-\eqref{war calkowe na mu} with
$$
c=0, \quad \nu_{G(0)}=0, \quad \mu(\dd v) =\mathbf{1}_{\{v>0\}}\frac{1}{v^{\alpha+1}}\dd v,  \quad \alpha\in(1,2).
$$
Our goal is to illustrate the following features of generating pairs $(G,Z)$ which do not appear in the case $d=1$:\\[1ex]
\noindent
$(a)$ for a given process $Z$ the function $G$ does not need to be unique, see Example \ref{ex1},\\[1ex]
\noindent
$(b)$ the coordinates of $Z$ may be of infinite variation, nevertheless,  $G(0)\neq 0$, see Example \ref{ex2}. 
In the case $d=1$ an important step in the proof of the form of \eqref{CIR equation geenralized}  was to show that $G(0)=0$, see Step 4 in the Proof of Theorem 2.1 in \cite{BarskiZabczykCIR}, see also Proposition 3.2 in \cite{BarskiZabczykArxiv}.\\[1ex]
\noindent
$(c)$ The coordinates of $Z$ do not need to be $\alpha$-stable, see Example \ref{ex3}.\\

We start with an $\alpha$-stable martingale in $\mathbb{R}^d, d>1$, with $\alpha\in(1,2)$ such that the spherical part $\lambda$ of its L\'evy measure is concentrated on 
$$
\mathbb{S}^{d-1}_{+}:= \{x\in\mathbb{R}^d:\mid x\mid=1, x\geq 0\}
$$
(writing $x \ge 0$ for $x \in \R^d$ we mean that all coordinates of $x$ are non-negative).
The Laplace exponent of the jump part of $Z$, identical with the Laplace exponent of $Z$,  admits the following representation:
\begin{align}\nonumber
J_\nu(z)&=\int_{\mathbb{S}^{d-1}_{+}}\lambda(\dd \xi)\int_{0}^{+\infty}\left(e^{-\langle z,r\xi\rangle}-1+\langle z,r\xi\rangle\right)\frac{1}{r^{1+\alpha}}\dd r\\[1ex]\nonumber
&=\int_{\mathbb{S}^{d-1}_{+}}\lambda(\dd \xi)\int_{0}^{+\infty}\left(e^{-r\langle z,\xi\rangle}-1+r\langle z,\xi\rangle\right)\frac{1}{r^{1+\alpha}}\dd r\\[1ex]\label{J dla alfa stabilnego}
&=C_{\alpha}\int_{\mathbb{S}^{d-1}_{+}}\langle z,\xi\rangle^\alpha\lambda(\dd \xi),
\end{align}
where $C_{\alpha}:=\Gamma(2-\alpha)/(\alpha(\alpha-1))$ and $\Gamma$ stands for the Gamma function.
Above we used the well known formula
$$
\int_{0}^{+\infty}\Big(e^{-uy}-1+uy\Big)\frac{1}{y^{1+\alpha}}\dd y=C_{\alpha} u^\alpha.
$$
\noindent
The following example shows that  the function $G$ in a generating pair $(G,Z)$ is not unique. 
\begin{ex}\label{ex1}
Let  $Z$ be an $\alpha$-stable martingale in $\mathbb{R}^d$ with the
Laplace exponent \eqref{J dla alfa stabilnego} and $G: [0,+\ns) \ra [0,+\ns)^d$, $G(0)=0$. Then a pair $(Z,G)$  generates an affine model if and only if the function $G$ satisfies
\begin{gather}\label{war na G alfa stabilny}
\int_{\mathbb{S}^{d-1}_{+}}\langle G(x),\xi\rangle^\alpha \lambda(\dd \xi)=\frac{C}{C_{\alpha}}x,\quad x\geq 0,
\end{gather}
with $C\ge 0$. 
We need to show that \eqref{war na G alfa stabilny} is equivalent to 
\eqref{war na exp Laplaca} with some measure $\mu(\dd v)$. Since $Z$ has no Wiener part and 
$\nu_{G(0)}(\dd v)\equiv 0$, we see that \eqref{war na exp Laplaca} takes the form
$$
J_{Z}(bG(x))=J_{\nu}(bG(x))=xJ_{\mu}(b), \quad x,b\geq 0.
$$
By \eqref{J dla alfa stabilnego} 
$$
J_\nu(bG(x))=C_{\alpha}\int_{\mathbb{S}^{d-1}_{+}}\langle bG(x),\xi\rangle^\alpha\lambda(\dd \xi)
=C_{\alpha} b^{\alpha}\int_{\mathbb{S}^{d-1}_{+}}\langle G(x),\xi\rangle^\alpha\lambda(\dd \xi).
$$
Consequently,
$$
C_{\alpha} b^{\alpha}\int_{\mathbb{S}^{d-1}_{+}}\langle G(x),\xi\rangle^\alpha\lambda(\dd \xi)=xJ_{\mu}(b),
$$
holds if and only if
$$
J_\mu(b)=C b^{\alpha}, \quad \int_{\mathbb{S}^{d-1}_{+}}\langle G(x),\xi\rangle^\alpha\lambda(\dd \xi)=\frac{C}{C_{\alpha}}x,
$$
for some $C \ge 0$. Hence, $\mu$ is an $\alpha$-stable measure and $G$ can be any function satisfying \eqref{war na G alfa stabilny}. 
\end{ex}

\begin{ex}\label{ex2}
Let $Z^{\alpha}(t)$ be a real valued $\alpha$-stable process with positive jumps only, $\alpha\in(1,2)$ and $G^{\alpha}(x):=x^{{1}/{\alpha}}$.
For 
$$
Z(t):=\Big(Z^{\alpha}(t),-Z^{\alpha}(t)\Big), \quad t\geq 0,
$$
$$
G(x):=\Big(G^\alpha(x)+1,-G^\alpha(x)+1\Big),\quad x\geq 0,
$$
and $F(x)\equiv 0$ equation \eqref{rownanie 2} becomes
\begin{align}\label{przyklad}\nonumber
\dd R(t)&=F(R(t))\dd t+\langle G(R(t-)),\dd Z(t)\rangle\\[1ex]\nonumber
&=\Big\langle(G^\alpha(R(t-))+1,-G^\alpha(R(t-))+1),(1,-1)\Big\rangle \dd Z^{\alpha}(t)\\[1ex]
&=2R(t-)^{{1}/{\alpha}}\dd Z^{\alpha}(t).
\end{align}
By \eqref{CIR equation geenralized} we see that $(G,Z)$ is a generating pair. Although the coordinates of $Z$ are of infinite variation,  $G(0)=(1,1)$. 

To see that $\nu_{G(0)}(\dd v)\equiv 0$ note that the L\'evy measure of $Z$ 
is supported by the half-line $\{ t(1,-1), t>0\}$  and therefore
$$
\langle G(0),y\rangle=\langle(1,1),(y_1,-y_1)\rangle=0, \quad y\in\emph{supp} \ \nu.
$$
It follows that
\begin{align*}
\nu_{G(0)}(A)&=\nu\{y\in \mathbb{R}^2: \langle G(0),y\rangle\in A\}\\[1ex]
&=\nu\{y\in \mathbb{R}^2: y_1+y_2\in A\} =0,
\end{align*}
provided that $0\notin A$.
\end{ex}
 Finally, we show that $Z$ does not need to have stable components.
\begin{ex}\label{ex3}
Let $E$ be any Borel subset of $[0, +\ns)$ such that 
$$
|E| = \int_{E} \dd r >0, \quad  \text{and} \quad |[0, +\ns) \setminus E| = \int_{[0, +\ns) \setminus E} \dd r  >0,
$$ 
and $Z_1$ and $Z_2$ be two independent L\'{e}vy processes with the L\'{e}vy measures 
\[
\nu_1(\dd r) = \mathbf{1}_E(r) \frac{\dd r}{r^{\alpha +1}}, \quad \nu_2(\dd r) = \mathbf{1}_{[0, +\ns) \setminus E}(r) \frac{\dd r}{r^{\alpha +1}}, \quad \alpha \in (1,2).
\]
Clearly, neither $Z_1$ nor $Z_2$ is stable, but $Z_1 + Z_2$ is, and has only positive jumps. Thus taking  $G(x) = \rbr{G_1(x), G_2(x)}= x^{1/\alpha} (1,1)$, $Z = \rbr{Z_1, Z_2}$ we get that the equation 
\[
\dd R(t) = \langle G(R(t-)),\dd Z(t)\rangle = R(t-)^{1/\alpha} \dd \rbr{Z_1(t)+Z_2(t)}
\]
generates an affine model. 
\end{ex}
Another example of a generating pair with $Z$ of independent and non-stable coordinates is presented in Section \ref{section Example in higher dimensions}.

\subsection{The form of affine models}

The following result has a supplementary character and shows how the functions $A(\cdot), B(\cdot)$ of the affine model \eqref{affine model}
are determined by a triplet  $(c,\nu_{G(0)}(\dd v),\mu(\dd v))$ satisfying  \eqref{nu G0 finite variation}-\eqref{war calkowe na mu}.

\begin{prop}\label{prop2}
Let the equation \eqref{rownanie 1} generate an affine model \eqref{affine model} with twice continuously differentiable functions $A(\cdot), B(\cdot)$. Let the drift $F(\cdot)$ and the projections $Z^{G(x)}$ satisfy 
\eqref{mult. CIR condition}, \eqref{rozklad nu G(x)} and \eqref{linear drift} with some constants $a,b,c$ and
measures $\nu_{G(0)}(dv)$, $\mu(dv)$. Then the functions $A,B$ are solutions of the following differential equations 
\begin{gather}\label{rr na B}
B^\prime(v)=aB(v)-\frac{1}{2}cB^2(v)-J_\mu(B(v))+1, \quad v\geq 0,\quad B(0)=0,
\end{gather}
\begin{gather}\label{A wzor od B}
A^\prime(v)=bB(v)-J_{\nu_{G(0)}}(B(v)), \quad v\geq 0, \quad A(0)=0.
\end{gather}
\end{prop}
\noindent
The proof of Proposition \ref{prop2} is postponed to Appendix.

\section{Noise with independent coordinates}\label{section Noise with independent coordinates}

This section deals with equation \eqref{rownanie 2} in the case when the coordinates $(Z_1,Z_2,...,Z_d)$ of the martingale $Z$ are independent processes. In view of Proposition \ref{prop wstepny}  we are interested in characterizing possible distributions of projections $Z^G$ over all generating pairs $(G,Z)$.  By \eqref{Z^G positive jumps} the jumps of the projections are necessarily positive. As the coordinates of $Z$ are independent, they do not jump together. Consequently, we see that, for each $x\geq 0$,
$$
\triangle Z^{G(x)}(t)=\langle G(x),\triangle Z(t)\rangle >0
$$
holds if and only if, for some $i=1,2,...,d$,
\begin{gather}\label{mucha}
G_i(x)\triangle Z_i(t)>0, \quad \triangle Z_j(t)=0, j\neq i.
\end{gather}
Condition \eqref{mucha} means that $G_i(x)$ and $\triangle Z_i(t)$ are of the same sign. We can consider only the case when both are positive, i.e.
$$
G_i(x)\geq 0, \quad i=1,2,...,d, \ x\geq 0, \qquad \triangle Z_i(t)\geq 0, \quad t> 0,
$$
because the opposite case can be turned into this one by replacing $(G_i,Z_i)$ with $(-G_i,-Z_i)$, $i=1,...,d$. The L\'evy measure $\nu_i(\dd y)$ of $Z_i$ is thus concentrated on $(0,+\infty)$ and, in view of \eqref{LaplaceZ}, the Laplace exponent of $Z_i$ takes the form
\begin{gather}\label{Laplace Zi}
J_i(b):=\frac{1}{2}q_{ii} b^2+\int_{0}^{+\infty}(e^{-b v}-1+b v)\nu_i(\dd v),\quad b \geq 0, \ i=1,2,...,d,
\end{gather}
with $q_{ii}\geq 0$. Recall, $q_{ii}$ stands on the diagonal of $Q$ - the covariance matrix of the Wiener part of $Z$.  
We will assume that $J_i, i=1,2,...,d$ are  {\it regularly varying at zero}. Recall, that means that
$$
\lim_{x\rightarrow 0^+}\frac{J_i(bx)}{J_i(x)}=\psi_i(b), \quad b> 0,\qquad i=1,2,...,d,
$$
for some function $\psi_i$. In fact $\psi_i$ is a power function, i.e.
$$
\psi_i(b)=b^{\alpha_i}, \quad b>0,
$$
with some $-\infty< \alpha_i<+\infty$ and $J_i$ is called to vary regularly with index $\alpha_i$. A characterization of slowly varying Laplace exponent in terms of the corresponding L\'evy measure is presented in Section \ref{sec slowly varying Laplace exp.}.

\subsection{Main results}

The main result of this section is the following.

\begin{tw} \label{TwNiez} Let $Z_1,...,Z_d$ be independent components of the L\'evy martingale $Z$ in $\R^d$. Assume that $Z_1,...,Z_{d}$  satisfy
\begin{equation} \label{ass1}
\triangle Z_i(t)\geq 0, \quad t>0, \quad Z_i \ \text{is of infinite variation}
\end{equation}
or 
\begin{equation} \label{ass2}
\triangle Z_i(t)\geq 0, \quad t>0, \text{ and } G(0)=0.
\end{equation}
Further, let us assume that for all $i=1,\ldots, d$ the Laplace exponent  \eqref{Laplace Zi} of $Z_i$ is regularly varying at $0$ and components of the function  $G$ satisfiy
$$
G_i(x)\geq 0, \ x\in[0,+\infty), \quad G_i \ \text{is continuous on } [0,+\infty).
$$
Then \eqref{rownanie 2} generates an affine model if and only 
 if $F(x)=ax+b$, $a\in\mathbb{R}, b\geq 0$, and the Laplace exponent $J_{Z^{G(x)}}$ of $Z^{G(x)}=\langle G(x), Z\rangle$ is of the form 
\begin{gather}\label{postac J_ZG przy niezaleznych}
J_{Z^{G(x)}}(b) = x \sum_{k=1}^g\eta_{k}b^{\alpha_{{k}}}, \quad  \eta_{k}> 0, \quad \alpha_k\in(1,2],  \quad k=1,2,\ldots,g,
\end{gather}
with some $1 \le g \le d$ and $\alpha_k\neq\alpha_j$ for $k\neq j$.
\end{tw}

Theorem \ref{TwNiez} allows determining the form of the measure $\mu(\dd v)$ in Proposition \ref{prop wstepny}.

\begin{cor}\label{cor o postaci mu} 
Let the assumptions of Theorem \ref{TwNiez} be satisfied. If equation \eqref{rownanie 2} generates an affine model 
then the function $J_\mu$ defined in \eqref{def Jmu i Jnu0} takes the form
\begin{gather}\label{J mu postaaaccccc}
J_{\mu}(b)  = \sum_{k=l}^{g}\eta_{k}b^{\alpha_{{k}}}, \quad l\in \{1,2\}, \quad  \eta_{k}> 0, \quad \alpha_k\in(1,2),  \quad k=l,l+1,\ldots,g,
\end{gather}
with $1 \le g \le d$, $\alpha_k\neq \alpha_j, k\neq j$ (for the case $l=2, g=1$ we set $J_{\mu}\equiv 0$, which means that $\mu(\dd v)$ disappears).

\end{cor}

Theorem \ref{TwNiez} specifies distributions of the projections $Z^{G(x)}$ of a generating pair $(Z,G)$.  As shown in Example \ref{ex3}, a given projection may correspond to many generating pairs  $(G,Z)$. This issue is also illustrated in Section \ref{section Generalized CIR equations on a plane} below, where all generating equations in the case $d=2$ are described. Below we show a tractable generating equation with the law of $Z^{G(x)}$ required by Theorem \ref{TwNiez}.

\begin{cor}
Let $R$ be the solution of \eqref{rownanie 2} with $F,G,Z$ satisfying the assumptions of Theorem \ref{TwNiez}. 
 Let $\tilde{Z}:=(\tilde{Z}_1,\tilde{Z}_2,...,\tilde{Z}_g)$ be a L\'evy martingale with independent stable coordinates
with indices $\alpha_k, k=1,2,...,g$, respectively, and $\tilde{G}(x)=(d_1 x^{1/\alpha_1},...,d_g x^{1/\alpha_g})$. 
Then 
$$
J_{Z^{G(x)}}(b)=J_{\tilde{Z}^{\tilde{G}(x)}}(b), \quad b,x\geq 0.
$$
Consequently, if $\tilde{R}$ is the solution of the equation 
\begin{gather}\label{rownanie sklajane}
\dd \tilde{R}(t)=(a\tilde{R}(t)+b) \dd t+\sum_{k=1}^{g}d_k^{1/{\alpha_k}} \tilde{R}(t-)^{1/{\alpha_k}}\dd \tilde{Z}_k(t),
\end{gather}
where $d_k:=({\eta_k}/{c_{k}})^{1/\alpha_k}, c_{k}=\frac{\Gamma(2-\alpha_k)}{\alpha_k(\alpha_k-1)}, k=1,2,...,g$, then the generators of $R$ and $\tilde{R}$ are equal.
\end{cor}
{\bf Proof:} By \eqref{postac J_ZG przy niezaleznych} we need to show that 
\begin{gather*}\label{rozklad rzutu konstrukcja}
J_{\tilde{Z}^{\tilde{G}(x)}}(b)=x\sum_{k=1}^{g}\eta_k b^{\alpha_k}, \quad b,x\geq 0.
\end{gather*}
Recall, the Laplace exponent of $\tilde{Z}_k$ equals
$J_k(b)=c_{k}b^{\alpha_k}, k=1,2,...,g$. By independence and the form of $\tilde{G}$ we have
\begin{align*}
J_{\tilde{Z}^{\tilde{G}(x)}}(b)&=\sum_{k=1}^{g}J_k(b\tilde{G}_k(x))=\sum_{k=1}^{g}c_k b^{\alpha_k}\frac{\eta_k}{c_k}x=x\sum_{k=1}^{g}\eta_k b^{\alpha_k}, \quad b,x\geq 0,
\end{align*}
as required. The second part of the thesis follows from Proposition \ref{prop wstepny}(B).\hfill$\square$

\subsubsection{Proofs}
The proofs of Theorem \ref{TwNiez} and Corollary \ref{cor o postaci mu} are preceded by two auxiliary results, i.e. Proposition \ref{bounds_alpha} and
Proposition \ref{rem o G(0)=0}. The first one provides some useful estimation for the function 
\begin{gather}\label{J}
J_{\rho}(b):=\int_{0}^{+\infty}(e^{-bv}-1+bv)\rho(\dd v), \quad b\geq 0,
\end{gather}
where the measure $\rho(\dd v)$ on $(0,+\ns)$ satisfies
\begin{gather}\label{nuJ}
0 < \int_0^{+\ns} \rbr{v^2\wedge v} \rho\rbr{\dd v} < +\ns.
\end{gather}
The second result shows that if all components of $Z$ are of infinite variation then $G(0)=0$.

\begin{prop} \label{bounds_alpha}
Let $J_{\rho}$ be a function given by \eqref{J} where the measure  $\rho$ satisfies \eqref{nuJ}. Then the function
$
(0,+\ns) \ni b \mapsto {J_{\rho}(b)}/{b}$ is strictly increasing and $\lim_{b \ra 0+}J_{\rho}(b)/b = 0$, while the function $(0,+\ns) \ni b \mapsto {J_{\rho}(b)}/{b^2}$
is strictly decreasing and $\lim_{b \ra +\ns}J_{\rho}(b)/b^2 = 0$. This yields, in particular, that, for any $b_0 >0$, 
\begin{gather}\label{oszacowania dwustronne J}
\frac{J_{\rho}\rbr{b_0}}{b_0^2}b^2 < J_{\rho}(b) < \frac{J_{\rho}\rbr{b_0}}{b_0}b, \quad b\in \rbr{0, b_0}.
\end{gather}
\end{prop} 
\noindent
{\bf Proof:} Let us start from the observation that the function 
$$
\frac{(1-e^{-t})t}{e^{-t}-1+t}, \quad t\geq 0,
$$
is strictly decreasing, with limit $2$ at zero and $1$ at infinity. This implies
 \begin{equation} \label{oszH}
 (e^{-t}-1+t) < (1-e^{-t})t < 2 (e^{-t}-1+t), \quad t \in (0, +\ns),
 \end{equation}
and, consequently,
$$
\int_{0}^{+\infty}(e^{-bv}-1+bv)\rho(\dd v) < \int_{0}^{+\infty}(1-e^{-bv})bv\ \rho(\dd v) < 2\int_{0}^{+\infty}(e^{-bv}-1+bv)\rho(\dd v), \quad b >0.
$$
This means, however, that
$$
J_{\rho}(b) < bJ_{\rho}^\prime(b) < 2J_{\rho}(b), \quad b > 0.
$$
So, we have
$$
\frac{1}{b} < \frac{J_{\rho}^\prime(b)}{J_{\rho}(b)}=\frac{d}{db}\ln J_{\rho}(b) < \frac{2}{b}, \quad b>0,
$$
and integration over some interval $[b_1,b_2]$, where $b_2 > b_1>0$, yields
$$
\ln b_2 - \ln b_1  < \ln J_{\rho}\rbr{b_2}-\ln J_{\rho}\rbr{b_1} < 2 \ln b_2 - 2 \ln b_1
$$
which gives that 
$$
\frac{J_{\rho}\rbr{b_2}}{b_2} > \frac{J_{\rho}\rbr{b_1}}{b_1}, \quad \frac{J_{\rho}\rbr{b_2}}{b_2^2} < \frac{J_{\rho}\rbr{b_1}}{b_1^2}.
$$

To see that $\lim_{b \ra 0+} {J_{\rho}\rbr{b}}/{b} = 0$ it is sufficient to use de l'H\^opital's rule,  \eqref{nuJ} and dominated convergence
$$
\lim_{b \ra 0+} \frac{J_{\rho}\rbr{b}}{b} = \lim_{b \ra 0+} {J'_{\rho}\rbr{b}} = \lim_{b \ra 0+}  \int_{0}^{+\infty}(1-e^{-bv}) v\ \rho(\dd v) = 0.
$$

To see that $\lim_{b \ra +\ns} {J_{\rho}\rbr{b}}/{b^2} = 0$ we also use de l'H\^opital's rule,  \eqref{nuJ} and dominated convergence.
If $\int_0^{+\ns} v\ \rho\rbr{\dd v} < +\ns$, then we have
$$
\lim_{b \ra +\ns} \frac{J_{\rho}\rbr{b}}{b^2} = \lim_{b \ra +\ns} \frac{J_{\rho}'\rbr{b}}{2b} = \frac{\int_0^{+\ns} v \rho\rbr{\dd v}}{+\ns}= 0.
$$
If $\int_0^{+\ns} v\ \rho\rbr{\dd v} = +\ns$ then we apply de l'H\^opital's rule twice and obtain
$$
\lim_{b \ra +\ns} \frac{J_{\rho}\rbr{b}}{b^2} = \lim_{b \ra +\ns} \frac{J_{\rho}'\rbr{b}}{2b} = \lim_{b \ra +\ns} \frac{J_{\rho}''\rbr{b}}{2} = \frac{1}{2} \lim_{b \ra +\ns} \int_{0}^{+\infty}e^{-bv} v^2\ \rho(\dd v) = 0.
$$
\hfill $\square$

\begin{prop}\label{rem o G(0)=0}
If $(G,Z)$ is a generating pair and all components of $Z$ are of infinite variation then $G(0)=0$. 
\end{prop}
\noindent 
{\bf Proof:}  Let $(G,Z)$ be a generating pair. Since the components of $Z$ are independent, its characteristic triplet \eqref{chrakterystyki Z} is such that $Q=\{q_{i,j}\}$ is a diagonal matrix, i.e.
$$
q_{ii}\geq 0, \quad q_{i,j}=0, \qquad i\neq j, \quad i,j=1,2,...,d,
$$
and the support of $\nu(\dd y)$ is contained in the positive half-axes of $\mathbb{R}^d$, see \cite{Sato} p.67. 
On the $i^{th}$ positive  half-axis 
\begin{gather}\label{rozlozenie ny}
\nu(\dd y)=\nu_i(dy_i),\qquad y=(y_1,y_2,...,y_d),
\end{gather}
for $i=1,2,...,d$.
 The $i^{th}$ coordinate of $Z$ is of infinite variation if and only if its Laplace exponent \eqref{Laplace Zi} is such that $q_{ii}>0$ or
\begin{gather}\label{rerere}
\int_{0}^{1}y_i\nu_i(\dd y_i)=+\infty,
\end{gather}
see \cite[Lemma 2.12]{Kyprianou}. It follows from \eqref{mult. CIR condition} that
$$
\frac{1}{2}\langle QG(x), G(x)\rangle=\frac{1}{2}\sum_{j=1}^{d}q_{jj}G_j^2(x)=cx,
$$
so if $q_{ii}>0$ then $G_i(0)=0$. If it is not the case, using \eqref{rozlozenie ny} and  \eqref{nu G0 finite variation} we see that the integral
\begin{align*}
\int_{(0,+\infty)}v\nu_{G(0)}(\dd v) & =\int_{\mathbb{R}^d_{+}}\langle G(0),y\rangle \nu(\dd y) \\ 
&  =\sum_{j=1}^{d}\int_{(0,+\infty)}G_j(0)y_j \ \nu_j(\dd y_j) =\sum_{j=1}^{d}
G_j(0) \int_{(0,+\infty)}y_j \ \nu_j(\dd y_j),
\end{align*}
is finite, so if  \eqref{rerere} holds then $G_i(0)=0$.\hfill $\square$

 \vskip2ex

\noindent
{\bf Proof of Theorem \ref{TwNiez}:} By assumption \eqref{ass1} and Proposition \ref{rem o G(0)=0} or by assumption \eqref{ass2}  we have $G(0)=0$, so it follows from Remark \ref{rem warunki w eksp. Laplacea} that
\begin{equation} \label{eq:dwa_zero}
J_{Z^{G(x)}}(b) = J_{1}(bG_{1}(x))+J_2(bG_2(x))+...+J_d(bG_d(x))=x\tilde{J}_{\mu}(b), \quad b,x\geq 0,
\end{equation} 
where $\tilde{J}_{\mu}(b) = c b^2 + {J}_{\mu}(b)$, $c\ge 0$ and ${J}_{\mu}(b)$ is given by \eqref{def Jmu i Jnu0}.
This yields
\begin{equation}
\frac{J_{1}\rbr{b\cdot G_{1}(x)}}{J_{1}\rbr{G_{1}(x)}}\cdot\frac{J_{1}\rbr{G_{1}(x)}}{x}+\ldots+\frac{J_{d}\rbr{b\cdot G_{d}(x)}}{J_{d}\rbr{G_{d}(x)}}\cdot\frac{J_{d}\rbr{G_{d}(x)}}{x}=\tilde{J}_{\mu}(b),\label{eq:dwa}
\end{equation}
where in the case $G_i(x) = 0$ we set $\frac{J_{i}\rbr{b\cdot G_{i}(x)}}{J_{i}\rbr{G_{i}(x)}}\cdot\frac{J_{i}\rbr{ G_{i}(x)}}{x} = 0$. Without loss of generality we may assume that $J_{1}$, $J_{2}$,$\ldots$,$J_{d}$ are non-zero (thus positive for positive arguments).  
By assumption, $J_{i}$, $i=1,2,\ldots, d$ vary regularly at $0$ with some indices $\alpha_{i}$, $i=1,2,\ldots,d$, 
so for $b>0$
\begin{equation} \label{eq:trzy}
\lim_{y\ra0+}\frac{J_{i}\rbr{b\cdot y}}{J_{i}(y)}=b^{\alpha_{i}}.
\end{equation}
Assume that 
\[
\alpha_{1}=\ldots=\alpha_{i\rbr{1}}>\alpha_{i\rbr{1}+1}=\ldots=\alpha_{i\rbr{2}}>\ldots\ldots>\alpha_{i\rbr{g-1}+1}=\ldots=\alpha_{i\rbr{g}}=\alpha_{d},
\]
where $i(g) = d$. Let us denote $i_0 = 0$ and
\begin{equation} \label{limits}
\eta_{k}(x) :=\frac{J_{i\rbr{k-1}+1}\rbr{G_{i\rbr{k-1}+1}(x)}+\ldots+J_{i\rbr{k}}\rbr{G_{i\rbr{k}}(x)}}{x}, \quad k=1,2,\ldots, g.
\end{equation}
We can rewrite equation \eqref{eq:dwa} in the form
\begin{equation}
\sum_{k=1}^{g} \rbr{\sum_{i=i\rbr{k-1}+1}^{i\rbr{k}} \frac{J_{i}\rbr{b\cdot G_{i}(x)}}{J_{i}\rbr{G_{i}(x)}}\cdot\frac{J_{i}\rbr{ G_{i}(x)}}{x}}=\tilde{J}_{\mu}(b),\label{eq:trzyy}
\end{equation}
By passing to the limit as $x\ra0+$, from \eqref{eq:trzy} and \eqref{eq:trzyy} we get 
\begin{align}
b^{\alpha_{i\rbr{1}}}\rbr{ \lim_{x \ra 0+} \eta_{1}(x)} +\ldots+   b^{\alpha_{i\rbr{g}}} \rbr{\lim_{x \ra 0+} \eta_{g}(x)}  =\tilde{J}_{\mu}(b), \label{eq:trzyyy}
\end{align} 
thus
\begin{gather}\label{J mu tilde sum power}
\tilde{J}_{\mu}(b)=\sum_{k=1}^g \eta_{k}b^{\alpha_{i\rbr{k}}},
\end{gather}
providing that the limits $\eta_{k} := \lim_{x \ra 0+} \eta_{k}(x)$, $k=1,2,\ldots, g$, exist. 
Thus it remains to prove that for  $k=1,2,\ldots, g$ the limits $\lim_{x \ra 0+} \eta_{k}(x)$ indeed exist {and that $\alpha_{i(k)} \in (1,2]$.}

First we will prove that $\lim_{x \ra 0+} \eta_{g}(x)$ exists.
Assume, by contrary, that this is not true, so 
\begin{equation} 
\limsup_{x \ra 0+} \eta_{g}(x) - \liminf_{x \ra 0+} \eta_{g}(x) \ge \delta >0.
\label{sequencess}
\end{equation}
It follows from \eqref{eq:dwa_zero} that
\[
\frac{J_1(G_1(x))+J_2(G_2(x))+...+J_d(G_d(x))}{x} = \sum_{k=1}^g \eta_k(x)=\tilde{J}_{\mu}(1).
\]
Let now $b_0 \in (0,1)$ be small enough so that 
\begin{equation} \label{oszacowanie}
\tilde{J}_{\mu}(1) b_0^{\alpha_{i\rbr{g-1}} - \alpha_{i(g)}} < \frac{\delta}{6}.
\end{equation}
Let us set in \eqref{eq:trzyy} $b=b_0$ and then divide both sides of \eqref{eq:trzyy} by $b_0^{\alpha_{i(g)}}$. 
For  $x>0$ sufficiently close to $0$ we have 
\[
\eta_g(x) - \frac{\delta}{6} \le \frac{1}{b_0^{\alpha_{i(g)}}} \rbr{ \sum_{i=i\rbr{g-1}+1}^{i\rbr{g}} \frac{{J}_{i}\rbr{b_0\cdot G_{i}(x)}}{J_{i}\rbr{G_{i}(x)}}\cdot\frac{J_{i}\rbr{ G_{i}(x)}}{x}} \le \eta_g(x) + \frac{\delta}{6}
\]
and 
\begin{align*}
\frac{1}{b_0^{\alpha_{i(g)}}}  \sum_{k=1}^{g-1}  \rbr{ \sum_{i=i\rbr{k-1}+1}^{i\rbr{k}} \frac{J_{i}\rbr{b_0\cdot G_{i}(x)}}{J_{i}\rbr{G_{i}(x)}}\cdot\frac{J_{i}\rbr{ G_{i}(x)}}{x}}  \le \sum_{k=1}^{g-1}  2 b_0^{\alpha_{i(k)} - \alpha_{i(g)}} \eta_k(x) \\
\le 2 b_0^{\alpha_{i(g-1)} - \alpha_{i(g)}} \tilde{J}_{\mu}(1)
\end{align*}
thus from \eqref{eq:trzyy}, two last estimates and \eqref{oszacowanie}
\[
\eta_g(x) - \frac{\delta}{6} \le \frac{\tilde{J}_{\mu}(b_0)}{b_0^{\alpha_{i(g)}}} \le \eta_g(x) + \frac{\delta}{6}+ 2\tilde{J}_{\mu}(1) b_0^{\alpha_{i(g-1)} - \alpha_{i(g)}} < \eta_g(x) + \frac{\delta}{2}.
\]
But this contradicts \eqref{sequencess}
since we must have 
\[
 \limsup_{x \ra 0+} \eta_g(x) \le \frac{\tilde{J}_{\mu}(b_0)}{b_0^{\alpha_{i(g)}}} + \frac{\delta}{6}, \quad \liminf_{x \ra 0+} \eta_g(x) \ge \frac{\tilde{J}_{\mu}(b_0)}{b_0^{\alpha_{i(g)}}} - \frac{\delta}{2}.
\]

Having proved the existence of the limits $\lim_{x \ra 0+} \eta_{g}(x)$, ..., $\lim_{x \ra 0+} \eta_{g-m+1}(x)$ we can proceed similarly to prove the existence of the limit $\lim_{x \ra 0+} \eta_{g-m}(x)$. 
Assume that $\lim_{x \ra 0+} \eta_{g-m}(x)$ does not exist, so \begin{equation} 
\limsup_{x \ra 0+} \eta_{g-m}(x) - \liminf_{x \ra 0+} \eta_{g-m}(x) \ge \delta >0.
\label{sequencess1}
\end{equation}
Let $b_0 \in (0,1)$ be small enough so that 
\begin{equation} \label{oszacowanie1}
\tilde{J}_{\mu}(1) b_0^{\alpha_{i\rbr{g-m-1}} - \alpha_{i(g-m)}} < \frac{\delta}{8}.
\end{equation}
Let us set in \eqref{eq:trzyy} $b=b_0$ and then divide both sides of \eqref{eq:trzyy} by $b_0^{\alpha_{i(g-m)}}$. 
For  $x>0$ sufficiently close to $0$ we have 
\[
\eta_{g-m}(x) - \frac{\delta}{8} \le \frac{1}{b_0^{\alpha_{i(g-m)}}} \sum_{i=i\rbr{g-m-1}+1}^{i\rbr{g-m}} \frac{{J}_{i}\rbr{b_0\cdot G_{i}(x)}}{J_{i}\rbr{G_{i}(x)}}\cdot\frac{J_{i}\rbr{ G_{i}(x)}}{x} \le \eta_{g-m}(x) + \frac{\delta}{8},
\]
\begin{align*}
\frac{1}{b_0^{\alpha_{i(g-m)}}}  \sum_{k=1}^{g-m-1}  \rbr{ \sum_{i=i\rbr{k-1}+1}^{i\rbr{k}} \frac{J_{i}\rbr{b_0\cdot G_{i}(x)}}{J_{i}\rbr{G_{i}(x)}}\cdot\frac{J_{i}\rbr{ G_{i}(x)}}{x}}  \le \sum_{k=1}^{g-m-1}  2 b_0^{\alpha_{i(k)} - \alpha_{i(g-m)}} \eta_k(x) \\
\le 2 b_0^{\alpha_{i(g-m-1)} - \alpha_{i(g-m)}} \tilde{J}_{\mu}(1)
\end{align*}
and
\begin{align*}
\sum_{k=g-m+1}^{g}  \frac{b_0^{\alpha_{i(k)}} \eta_k}{b_0^{\alpha_{i(g-m)}}} - \frac{\delta}{8} & \le \frac{1}{b_0^{\alpha_{i(g-m)}}} \sum_{k=g-m+1}^{g}  \sum_{i=i\rbr{k-1}+1}^{i\rbr{k}} \frac{{J}_{i}\rbr{b_0\cdot G_{i}(x)}}{J_{i}\rbr{G_{i}(x)}}\cdot\frac{J_{i}\rbr{ G_{i}(x)}}{x} \\
& \le \sum_{k=g-m+1}^{g}  \frac{b_0^{\alpha_{i(k)}} \eta_k}{b_0^{\alpha_{i(g-m)}}}  + \frac{\delta}{8}
\end{align*}
thus from \eqref{eq:trzyy}, last three estimates and \eqref{oszacowanie1}
\begin{align*}
\eta_{g-m}(x) - \frac{\delta}{4} & \le \frac{J_{\mu}(b_0)}{b_0^{\alpha_{i(g-m)}}} - \sum_{k=g-m+1}^{g}  \frac{b_0^{\alpha_{i(k)}} \eta_k}{b_0^{\alpha_{i(g-m)}}} \\
 & \le \eta_{g-m}(x) + \frac{\delta}{4}+ 2\tilde{J}_{\mu}(1) b_0^{\alpha_{i(g-1)} - \alpha_{i(g)}} <\eta_{g-m}(x) +  \frac{\delta}{2}.
\end{align*}
But this contradicts \eqref{sequencess1}.

Now we are left with the proof that for $k=1,2,\ldots,g$, $\alpha_{i(k)} \in (1,2]$. Since the Laplace exponent of $Z_i$ is given by 
\eqref{Laplace Zi}, by Proposition \ref{bounds_alpha} we necessarily have that $J_i$ varies regularly with index $\alpha_i \in [1,2], i=1,2,...,d$. Thus it remains to prove that $\alpha_i >1, i=1,2,...,d$. If it was not true we would have $\alpha_{i(g)}=1$ in \eqref{J mu tilde sum power} and $\eta_g>0$. Then
$$
\lim_{b\ra 0+} \tilde{J}_{\mu}(b)/b = \lim_{b\ra 0+} J_{\mu}(b)/b =\eta_{g}>0,
$$ 
but, again, by Proposition \ref{bounds_alpha} it is not possible. 
\hfill $\square$

\vskip2ex
\noindent
{\bf Proof of Corollary \ref{cor o postaci mu} :} From Remark \ref{rem warunki w eksp. Laplacea} and Theorem \ref{TwNiez} we know that 
\[
J_{Z^{G(x)}}(b) =  x  c b^2 + x{J}_{\mu}(b) = x \sum_{k=1}^g\eta_{k}b^{\alpha_{{k}}}, 
\]
where $1\le g \le d$, $\eta_{k}> 0$, $\alpha_k\in(1,2]$, $\alpha_k\neq\alpha_j$, $k,j=1,2,\ldots,g$, $c\ge 0$.
Without loss of generality we may assume that $2 \ge \alpha_1 > \alpha_2 > \ldots > \alpha_g >1$. Thus, since the Laplace exponent is nonnegative, $ x{J}_{\mu}(b)$ is of the form
\begin{gather}\label{pierwsza postac}
x{J}_{\mu}(b) = x\sum_{k=1}^g\eta_{k}b^{\alpha_{{k}}}, \qquad \text {if} \ c=0,
\end{gather}
or
\begin{gather}\label{druga postac}
x{J}_{\mu}(b) = x\sbr{(\eta_1-c)b^2+ \sum_{k=2}^g\eta_{k}b^{\alpha_{{k}}}}, \qquad \text {if} \ 0<c\leq\eta_1 \ \text {and} \ \alpha_1=2.
\end{gather}
In the case \eqref{pierwsza postac} we need to show that $\alpha_1<2$. If it was not true, we would have
\[
\lim_{b \ra +\ns} \frac{{J}_{\mu}(b)}{b^2} = \eta_1 >0,
\]
but this contradicts Proposition \ref{bounds_alpha}. In the same way we prove that $\eta_1=c$ in \eqref{druga postac}.
This proves the required representation \eqref{J mu postaaaccccc}.
\hfill $\square$

\subsection{Characterization of regularly varying Laplace exponents}\label{sec slowly varying Laplace exp.}

In this section we reformulate the assumption that $J_i, i=1,...,d$, vary regularly at zero 
in terms of the behaviour of the L\'evy measures of $Z_i, i=1,...,d$. As our considerations are componentwise,
we write for simplicity $\nu(\dd v):=\nu_i(\dd v)$ for the L\'evy measure of $Z_i$ and $J:=J_i$ for its Laplace exponent.

\begin{prop}\label{prop J ragularly varying warunki na nu}

Let $\nu(\dd v)$ be such that 
\begin{gather}\label{war calk w prop charakt.}
\int_{0}^{+\infty}(y^2\wedge y) \ \nu(dy)<+\infty.
\end{gather}
Let $\tilde{\nu}(\dd v)$ be the measure
$$
\tilde{\nu}(\dd v):=v^2\nu(\dd v),
$$
and $\tilde{F}$ its cumulative distribution function, i.e.
$$
\tilde{F}(v):=\tilde{\nu}((0,v))=\int_{0}^{v}u^2\nu(\dd u),\quad v\geq 0. 
$$
Then, for $\alpha\in(1,2)$, the following conditions are equivalent
\begin{gather}\label{J varying regggg.}
\lim_{x\rightarrow 0^{+}}\frac{J(bx)}{J(x)}=b^{\alpha}, \quad b\geq 0,
\end{gather}
$$
\lim_{y\rightarrow +\infty}\frac{\tilde{F}(by)}{\tilde{F}(y)}=b^{2-\alpha}, \quad b\geq 0.
$$
If, additionally, $\nu(\dd v)$ has a density function $g(v)$ such that 
\begin{gather}\label{niecalkowalnosc y2}
\int_{0}^{+\infty}v^2 g(v) \nu(\dd v)=+\infty,
\end{gather}
then \eqref{J varying regggg.} is equivalent 
to the condition
$$
\lim_{y\rightarrow +\infty}\frac{g(by)}{g(y)}=b^{-\alpha-1}, \quad b > 0.
$$
\end{prop}
{\bf Proof:}  Under \eqref{war calk w prop charakt.} the function $J$ given by \eqref{J} is well defined for $b\geq 0$, twice differentiable and 
\begin{gather*}
J^\prime(b)=\int_{0}^{+\infty}v(1-e^{-bv})\nu(\dd v), \quad 
J^{\prime\prime}(b)=\int_{0}^{+\infty}v^2e^{-bv}\nu(\dd v), \quad b\geq 0,
\end{gather*}
see \cite{Rusinek}, Lemma 8.1 and Lemma 8.2.
This implies that
\begin{align*}
\lim_{x\rightarrow 0^{+}}\frac{J(bx)}{J(x)}&=b\cdot \lim_{x\rightarrow 0^{+}}\frac{J^\prime(bx)}{J^\prime(x)}=
b^2\cdot \lim_{x\rightarrow 0^{+}}\frac{J^{\prime\prime}(bx)}{J^{\prime\prime}(x)}\\[1ex]
&=b^2\cdot \lim_{x\rightarrow 0^{+}}\frac{\int_{0}^{+\infty}e^{-bxv}v^2\nu(\dd v)}{\int_{0}^{+\infty}e^{-xv}v^2\nu(\dd v)}.
\end{align*}
Consequently, by \eqref{J varying regggg.} 
\begin{gather}\label{do TAubera}
\lim_{x\rightarrow 0^{+}}\frac{\int_{0}^{+\infty}e^{-bxv}v^2\nu(\dd v)}{\int_{0}^{+\infty}e^{-xv}v^2\nu(\dd v)}=b^{\alpha-2}.
\end{gather}
Notice, that the left side  is a quotient of two transforms of the measure $\tilde{\nu}(\dd v)$.
By the Tauberian theorem we have that \eqref{do TAubera} holds if and only if
$$
\frac{\tilde{F}(by)}{\tilde{F}(y)}\underset{y\rightarrow +\infty}{\longrightarrow}b^{2-\alpha}, \quad b\geq 0.
$$
If $\nu(\dd v)$ has a density $g(v)$ satisfying \eqref{niecalkowalnosc y2} then
\begin{align*}
\lim_{y\rightarrow +\infty}\frac{\tilde{F}(by)}{\tilde{F}(y)}&=\lim_{y\rightarrow +\infty}\frac{\int_{0}^{by}u^2g(u)\dd u}{\int_{0}^{y}u^2g(u)\dd u}
=\lim_{y\rightarrow +\infty}\frac{b\cdot (by)^2g(by)}{y^2g(y)}\\[1ex]
&=b^3\cdot\lim_{y\rightarrow +\infty}\frac{g(by)}{g(y)}.
\end{align*}
It follows that
$$
\lim_{y\rightarrow +\infty}\frac{g(by)}{g(y)}=b^{-\alpha-1}.
$$
which proves the result.\hfill$\square$

\begin{rem}
By general characterization of regularly varying functions we see that the functions
$\tilde{F}$ and $g$ from Proposition \ref{prop J ragularly varying warunki na nu} must be of the forms
$$
\tilde{F}(b)=b^{2-\alpha}L(b), \quad b\geq 0,
$$
$$
g(b)=b^{-\alpha-1}\tilde{L}(b), \quad b\geq 0,
$$
where $L$ and $\tilde{L}$ are slowly varying functions at $+\infty$, i.e.
$$
\frac{L(by)}{L(y)}\underset{y\rightarrow +\infty}{\longrightarrow}1, \quad \frac{\tilde{L}(by)}{\tilde{L}(y)}\underset{y\rightarrow +\infty}{\longrightarrow}1.
$$
\end{rem}

\subsection{Generating equations on a plane}\label{section Generalized CIR equations on a plane}

In this section we characterize all equations \eqref{rownanie 2}, with $d=2$, which generate affine models.
In view of Theorem \ref{TwNiez} generating pairs $(G,Z)$ are such that 
\begin{gather}\label{równanie z Filipovica d=2}
J_1(bG_1(x))+J_2(bG_2(x))=x \tilde{J}_{\mu}(b), \quad b,x\geq 0,
\end{gather}
where $\tilde{J}_{\mu}$ takes one of the two following forms
\begin{gather}\label{pierwszy przyp Jmu}
\tilde{J}_{\mu}(b)=\eta_1 b^{\alpha_1}, \quad b\geq 0,
\end{gather}
or
\begin{gather}\label{drugi przyp Jmu}
\tilde{J}_{\mu}(b)=\eta_1 b^{\alpha_1}+\eta_2 b^{\alpha_2},\quad b\geq 0,
\end{gather}
where $\eta_1,\eta_2>0$, $2\geq \alpha_1>\alpha_2>1$. We deduce from \eqref{równanie z Filipovica d=2} the form of $G$ and characterize the noise $Z$.

\begin{tw}\label{tw d=2 independent coord.}
Let $G(x)=(G_1(x), G_2(x))$ be continuous functions such that $G_1(x)>0,G_2(x)>0, x>0$ and $\frac{G_2(x)}{G_1(x)}\in C^1(0,+\infty)$.
Let $Z(t)=(Z_1(t),Z_2(t))$ have independent coordinates of infinite variation with Laplace exponents varying regularly
with indices $\alpha_1,\alpha_2$, respectively, where $2\geq \alpha_1\geq\alpha_2>1$. 
\begin{enumerate}[I)]
\item If $\tilde{J}_{\mu}$  is of the form \eqref{pierwszy przyp Jmu} 
then $(G,Z)$ is a  generating pair if and only if one of the following two cases holds:
\begin{enumerate}[a)]
\item 
\begin{gather}\label{wspolinionwosc G}
G(x)=c_0 \ x^{1/\alpha_1}\cdot\left(
\begin{array}{ccc}
 G_1\\
 G_2, 
\end{array}
\right), \quad x\geq 0,
\end{gather}
where $c_0>0, G_1>0, G_2>0$ and the process
$$
G_1 Z_1(t)+G_2 Z_2(t), \quad t\geq 0,
$$
is $\alpha_1$-stable.
\item $G(x)$ is such that
\begin{gather}\label{Ib}
c_1 G^{\alpha_1}_1(x)+c_2 G^{\alpha_1}_2(x)=\eta_1 x, \quad x\geq 0,
\end{gather}
with some constants $c_1,c_2>0$, and $Z_1,Z_2$ are $\alpha_1$-stable processes.
\end{enumerate}

\item If $\tilde{J}_{\mu}$  is of the form \eqref{drugi przyp Jmu}, then $(G,Z)$ is a  generating pair if and only if
\begin{gather}\label{postac g w drugim prrzypadku}
G_1(x)=\left(\frac{\eta_1}{c_1}  x\right)^{1 / \alpha_1}, \quad G_2(x)=\left(\frac{\eta_2}{c_2}  x\right)^{1 / \alpha_2}, \quad x\geq 0,
\end{gather}
with $c_1,c_2>0$ and $Z_1$ is $\alpha_1$-stable, $Z_2$ is $\alpha_2$-stable.
\end{enumerate}
\end{tw}

\noindent
{\bf Proof:} First let us consider the case when 
\begin{gather}\label{znikanie pochodnej}
\left(\frac{G_2(x)}{G_1(x)}\right)^\prime=0, \qquad x>0.
\end{gather}
Then 
$G(x)$ can be written in the form
\begin{gather*}
G(x)=g(x)\cdot\left(
\begin{array}{ccc}
 G_1\\
 G_2, 
\end{array}
\right), \quad x\geq 0,
\end{gather*}
with some function $g(x)\geq 0, x\geq 0$, and constants $G_1>0,G_2>0$. Equation \eqref{rownanie 2} 
amounts then to
\begin{align*}
dR(t)&=F(R(t))+g(R(t-)) \left(G_1dZ_1(t)+G_2 dZ_2(t)\right)\\[1ex]
&=F(R(t))+g(R(t-)) d\tilde{Z}(t), \quad t\geq 0,
\end{align*}
which is an equation driven by the one dimensional L\'evy process $\tilde{Z}(t):=G_1 Z_1(t)+G_2 Z_2(t)$. 
It follows that $\tilde{Z}$ is $\alpha_1$-stable with $\alpha_1\in(1,2]$ and that $g(x)=c_0x^{1/ \alpha_1}, c_0>0$. 
Notice that $Z^{G(x)}(t)=c_0 x^{\frac{1}{\alpha_1}}\tilde{Z}$, so $J_{Z^{G(x)}}(b)=C_{\alpha_1}(c_0 x^{\frac{1}{\alpha_1}}b)^{\alpha_1}=x c_0^{\alpha_1}C_{\alpha_1} b^{\alpha_1}$.
Hence \eqref{pierwszy przyp Jmu} holds and this proves $(Ia)$.

If \eqref{znikanie pochodnej} is not satisfied, then
\begin{gather}\label{niezerowanie pochodnej}
\left(\frac{G_2(x)}{G_1(x)}\right)^\prime\neq 0, \quad x\in (\underline{x},\bar{x}),
\end{gather}
in some interval $(\underline{x},\bar{x})\subset (0,+\infty)$. In the rest of the proof we consider this case and
prove $(Ib)$ and $(II)$.

$(Ib)$ From the equation
\begin{gather}\label{rrr}
J_1(bG_1(x))+J_2(bG_2(x))=x\eta_1 b^{\alpha_1}, \quad b\geq 0, \ x\geq 0,
\end{gather}
we explicitely determine unknown functions. Inserting $b/G_1(x)$ for $b$ yields
\begin{gather}\label{rowwww do eliminacji J1}
J_1(b)+J_2\left(b\frac{G_2(x)}{G_1(x)}\right)=\eta_1\frac{x}{G_1^{\alpha_1}(x)}b^{\alpha_1}, \quad b\geq 0, \quad x>0.
\end{gather}
Differentiation over $x$ yields
 $$
J_2^\prime\left(b\frac{G_2(x)}{G_1(x)}\right)\cdot b \left(\frac{G_2(x)}{G_1(x)}\right)^\prime
=\eta_1\left(\frac{x}{G_1^{\alpha_1}(x)}\right)^\prime b^{\alpha_1}, \quad b\geq 0 ,\quad  x>0.
 $$
Using \eqref{niezerowanie pochodnej} and dividing by $\left(\frac{G_2(x)}{G_1(x)}\right)^\prime$ leads to
 $$
J_2^\prime\left(b\frac{G_2(x)}{G_1(x)}\right)\cdot b=\eta_1 \frac{\left(\frac{x}{G^{\alpha_1}_1(x)}\right)^\prime}{\left(\frac{G_2(x)}{G_1(x)}\right)^\prime}\cdot b^{\alpha_1},\quad b\geq 0 ,\quad x\in(\underline{x},\bar{x}).
$$
By inserting $b\frac{G_1(x)}{G_2(x)}$ for $b$ one computes the derivative of $J_2$: 
$$
J_2^\prime(b)=\eta_1 \frac{\left(\frac{x}{G^{\alpha_1}_1(x)}\right)^\prime\left(\frac{G_1(x)}{G_2(x)}\right)^{\alpha_1-1}}{\left(\frac{G_2(x)}{G_1(x)}\right)^\prime}\cdot b^{\alpha_1-1},\quad b>0 , \quad x\in(\underline{x},\bar{x}).
$$
Fixing $x$ and integrating over $b$ provides
\begin{gather}\label{J_2 wylioczona}
J_2(b)=c_2 b^{\alpha_1}, \quad b>0,
\end{gather}
with some $c_2\geq 0$. Actually $c_2>0$ as $Z_2$ is of infinite variation and $J_2$ can not disappear.

By the symmetry of \eqref{rrr} the same conclusion holds for $J_1$, i.e.
\begin{gather}\label{J_1 wylioczona}
J_1(b)=c_1 b^{\alpha_1}, \quad b>0,
\end{gather}
with $c_1>0$. Using \eqref{J_2 wylioczona} and \eqref{J_1 wylioczona} in \eqref{rrr} gives us \eqref{Ib}. This proves $(Ib)$.

 $II)$ Solving the equation 
\begin{gather}\label{rrrrr}
J_1(bG_1(x))+J_2(bG_2(x))=x(\eta_1 b^{\alpha_1}+\eta_2 b^{\alpha_2}), \quad b,x\geq 0,
\end{gather}
in the same way as we solved \eqref{rrr} yields that 
 \begin{gather}\label{J1,J2 podwojne stabilne}
 J_1(b)=c_1 b^{\alpha_1}+c_2 b^{\alpha_2}, \quad  J_2(b)=d_1 b^{\alpha_1}+d_2 b^{\alpha_2}, \quad b\geq 0,
 \end{gather}
 with $c_1,c_2,d_1,d_2\geq 0$, $c_1+c_2>0, d_1+d_2>0$. From \eqref{rrrrr} and \eqref{J1,J2 podwojne stabilne} we can specify the following conditions for $G$:
  \begin{align}\label{aaaaaa}
c_1G_1^{\alpha_1}(x)+d_1G_2^{\alpha_1}(x)&=\eta_1 x,\\[1ex]\label{bbbbbb}
c_2G_1^{\alpha_2}(x)+d_2G_2^{\alpha_2}(x)&=\eta_2 x.
\end{align}
 We show that $c_1>0, c_2=0, d_1=0, d_2>0$ by excluding the opposite cases.

 If $c_1>0,c_2>0$, one computes from \eqref{aaaaaa}-\eqref{bbbbbb} that
\begin{gather}\label{G2 na dwa sposoby}
G_1(x)=\left(\frac{1}{c_1}(\eta_1x-d_1G_2^{\alpha_1}(x))\right)^{\frac{1}{\alpha_1}}=\left(\frac{1}{c_2}(\eta_2x-d_2G_2^{\alpha_2}(x))\right)^{\frac{1}{\alpha_2}}, \quad x\geq 0.
\end{gather}
 This means that, for each $x\geq 0$, the value $G_2(x)$ is a solution of the following equation in the $y$-variable
\begin{gather}\label{rownanie z y}
\left(\frac{1}{c_1}(\eta_1x-d_1y^{\alpha_1})\right)^{\frac{1}{\alpha_1}}=\left(\frac{1}{c_2}(\eta_2x-d_2y^{\alpha_2})\right)^{\frac{1}{\alpha_2}}, 
\end{gather} 
with $y\in \left[0,\left(\frac{\gamma_1 x}{d_1}\right)^{\frac{1}{\alpha_1}}\wedge
\left(\frac{\gamma_2 x}{d_2}\right)^{\frac{1}{\alpha_2}}\right]$.
If $d_1=0$ or $d_2=0$ we compute $y=y(x)$ from \eqref{rownanie z y} and see that its positivity is broken close to zero  or for large $x$. We need to exclude the case $d_1>0,d_2>0$. However,  in the case $c_1, c_2,d_1,d_2>0$ equation \eqref{rownanie z y} has no solution because, for large $x>0$, the left side of \eqref{rownanie z y} is strictly less then the right side. This inequality follows from Proposition \ref{prop o braku rozwiazan} proven below. 

So, we proved that $c_1\cdot c_2=0$ and similarly one proves that $d_1\cdot d_2=0$.
The case $c_1=0, c_2>0, d_1>0, d_2=0$ can be rejected because then $J_1$ would vary regularly with index $\alpha_2$
and $J_2$ with index $\alpha_1$, which is a contradiction. It follows that $c_1>0, c_2=0, d_1=0, d_2>0$ and in this case we obtain \eqref{postac g w drugim prrzypadku} from \eqref{aaaaaa} and \eqref{bbbbbb}. \hfill$\square$
\begin{prop}\label{prop o braku rozwiazan}
Let  $a,b,c,d>0$, $\gamma\in(0,1)$,  $2\geq \alpha_1>\alpha_2>1$. Then for large $x>0$ the following inequalities are true
\begin{gather}\label{niernier pomocnicza}
\Big(ax-(bx-cz)^\gamma\Big)^{\frac{1}{\gamma}}-dz> 0, \qquad z\in \Big[0,\frac{b}{c}x\Big],
\end{gather}
\begin{gather}\label{niernier pomocnicza wlasciwa}
(bx-cy^{\alpha_1})^{\frac{1}{\alpha_1}}<(ax-dy^{\alpha_2})^{\frac{1}{\alpha_2}}, \quad y\in \Big[0,\Big(\frac{b}{c}x\Big)^{\frac{1}{\alpha_1}}\wedge \Big(\frac{a}{d}x\Big)^{\frac{1}{\alpha_2}}\Big].
\end{gather}
\end{prop}
{\bf Proof:} First we prove \eqref{niernier pomocnicza} and write it in the equivalent form
\begin{gather}\label{equivalentt}
ax\geq (dz)^{\gamma}+(bx-cz)^\gamma=:h(z).
\end{gather}
Since
$$
h^{\prime}(z)=\gamma \Big(d^\gamma z^{\gamma-1}-c(bx-cz)^{\gamma-1}\Big),
$$
$$
h^{\prime\prime}(z)=\gamma(\gamma-1)\Big(d^\gamma z^{\gamma-2}+c^2(bx-cz)^{\gamma-2}\Big)< 0, \quad z\in \Big[0,\frac{b}{c}x\Big],
$$
the function $h$ is concave and attains its maximum at point
$$
z_0:=\theta x:=\frac{b c^{\frac{1}{\gamma-1}}}{d^{\frac{\gamma}{\gamma-1}}+c^{\frac{\gamma}{\gamma-1}}}x \in \Big[0,\frac{b}{c}x\Big],
$$
which is a root of $h^\prime$. It follows that
\begin{align*}
h(z)\leq h(\theta x)&=(\theta x)^{\gamma}+(bx-c\theta x)^{\gamma}\\
&=(\theta^\gamma+(b-c\theta)^\gamma)x^{\gamma}< ax.
\end{align*}
The last strict inequality holds for large $x$ and \eqref{niernier pomocnicza} follows. \eqref{niernier pomocnicza wlasciwa} follows from
\eqref{niernier pomocnicza} by setting $\gamma=\alpha_2/\alpha_1$, $z=y^{\alpha_1}$.
\hfill$\square$

\subsection{An example in higher dimensions}\label{section Example in higher dimensions}

In Theorem \ref{tw d=2 independent coord.} we showed that in the case $d=2$ there exists a unique 
generating equation corresponding to the function 
\begin{gather}\label{drugi przyp Jmu 2}
\tilde{J}_{\mu}(b)=\eta_1 b^{\alpha_1}+\eta_2 b^{\alpha_2},\quad b\geq 0,
\end{gather}
with $\eta_1,\eta_2>0, 2\geq \alpha_1>\alpha_2>1$. This function does not guarantee the uniqueness of generating equations in higher dimensions. Below we show a family of generating pairs $(G,Z)$ taking values in $\mathbb{R}^3$ such that $J_{Z^{G(x)}}(b)=x\tilde{J}_{\mu}(b)$.

\begin{ex} 
Let us consider a process $Z(t)=(Z_1(t),Z_2(t),Z_3(t))$ with independent coordinates such that  $Z_1$ is $\alpha_1$-stable, $Z_2$ is $\alpha_2$-stable,  $Z_3$ is a sum of an $\alpha_1$-  and $\alpha_2$-stable processes. Then
$$
J_1(b)=\gamma_1 b^{\alpha_1}, \quad J_2(b)=\gamma_2 b^{\alpha_2}, \quad J_3(b)=\gamma_3b^{\alpha_1}+\tilde{\gamma}_3 b^{\alpha_2}, \quad b\geq 0,
$$
where $\gamma_1>0,\gamma_2>0,\gamma_3>0,\tilde{\gamma}_3>0$. We are looking for non-negative functions $G_1, G_2,G_3$ solving the equation
\begin{gather}\label{krokodyl}
J_1(bG_1(x))+J_2(bG_2(x))+J_3(bG_3(x))=x \tilde{J}_\mu(b), \quad x,b\geq 0,
\end{gather}
where $\tilde{J}_{\mu}$ is given by \eqref{drugi przyp Jmu 2}.
It follows from \eqref{krokodyl} that
$$
\gamma_1b^{\alpha_1}(G_1(x))^{\alpha_1}+\gamma_2b^{\alpha_2}(G_2(x))^{\alpha_2}
+\gamma_3b^{\alpha_1}(G_3(x))^{\alpha_1}+\tilde{\gamma}_3 b^{\alpha_2}(G_3(x))^{\alpha_2}=
x\left[\eta_1 b^{\alpha_1}+\eta_2b^{\alpha_2}\right], \quad x,b\geq 0,
$$
and, consequently,
$$
b^{\alpha_1}\left[\gamma_1G_1^{\alpha_1}(x)+\gamma_3 G_3^{\alpha_1}(x)\right]+
b^{\alpha_2}\left[\gamma_2G_2^{\alpha_2}(x)+\tilde{\gamma}_3 G_3^{\alpha_2}(x)\right]=x\left[\eta_1 b^{\alpha_1}+\eta_2b^{\alpha_2}\right], \quad x,b\geq 0.
$$
Thus we obtain the following system of equations 
\begin{gather*}
\gamma_1G_1^{\alpha_1}(x)+\gamma_3 G_3^{\alpha_1}(x)=x\eta_1, \\
\gamma_2G_2^{\alpha_2}(x)+\tilde{\gamma}_3 G_3^{\alpha_2}(x)=x\eta_2,
\end{gather*}
which allows us to determine $G_1$ and $G_2$ in terms of $G_3$, that is
\begin{gather}\label{G1}
G_1(x)=\left(\frac{1}{\gamma_1}\left(x\eta_1-\gamma_3 G_3^{\alpha_1}(x)\right)\right)^{\frac{1}{\alpha_1}}\\\label{G2}
G_2(x)=\left(\frac{1}{\gamma_2}\left(x\eta_2-\tilde{\gamma}_3 G_3^{\alpha_2}(x)\right)\right)^{\frac{1}{\alpha_2}}.
\end{gather}
The positivity of $G_1,G_2,G_3$  means that $G_3$ satisfies
\begin{gather}\label{war nieuj}
0\leq G_3(x)\leq \left(\frac{\eta_1}{\gamma_3}x\right)^{\frac{1}{\alpha_1}}\wedge \left(\frac{\eta_2}{\tilde{\gamma}_3}x\right)^{\frac{1}{\alpha_2}}, \quad x\geq 0.
\end{gather}
It follows that $(G,Z)$ with any $G_3$ satisfying \eqref{war nieuj} and $G_1,G_2$ given by \eqref{G1}, \eqref{G2} constitutes a generating pair.
\end{ex}

\section{Spherical L\'evy noise}\label{section Spherical Levy noise}

This section deals with equation \eqref{rownanie 2} in the case when $Z$ is a spherical L\'evy process with characteristic triplet
$(a,Q,\nu(\dd y))$. Recall, the L\'evy measure $\nu(\dd y)$ is described in terms of a finite spherical measure $\lambda(\dd \xi)$ on $\mathbb{S}^{d-1}$ and a radial measure $\gamma(\dd r)$ on $(0,+\infty)$ by the formula \eqref{spherical}. 
As $Z$ is a martingale, by \eqref{chrakterystyki Z}, 
$$
a=-\int_{\mid y\mid>1}\mid y\mid\nu(dy)= \int_{\S^{d-1}}\lambda(\dd \xi)\int_{1}^{+\infty}r\xi \ \gamma(\dd r),
$$
and the integrability of $Z$ implies that
\begin{gather}\label{warunek na miare radialna a}
\int_{\mid y\mid>1}\mid y\mid\nu(dy)=\int_{\S^{d-1}}\lambda(\dd \xi)\int_{0}^{+\infty}\mid r\xi\mid\gamma(dr)=\lambda(\S^{d-1})\cdot \int_{1}^{+\infty}r \gamma(dr) <+\infty.
\end{gather}
The jump part of $Z$ is assumed to have infinite variation, which means that 
$$
\int_{\mid y\mid\leq 1}\mid y\mid\nu(dy)=\lambda(\S^{d-1})\cdot \int_{0}^{1}r \gamma(dr)=+\infty.
 $$
Consequently,  the radial measure is of infinite variation, i.e.
\begin{gather}\label{warunek na miare radialna b}
\int_{0}^{1}r \gamma(dr)=+\infty.
\end{gather}

Furthermore, the measure $\nu(\dd y)$ will be assumed to have its support not contained in any proper linear subspace of $\R^d$, i.e.
\begin{gather}\label{spherical part non-degenerated}
\text{Linear span  (supp } \lambda) = \R^d.
\end{gather}
Moreover, by Proposition \ref{prop wstepny}(A(a)) and \eqref{spherical},
$
\lambda \cbr{\xi \in \S^{d-1}: \langle G(0), \xi \rangle <0} = 0
$
which implies that 
\begin{equation} \label{nonzerol}
\langle G(0),\xi \rangle \ge 0 \ \text{for any } \xi \in \text{supp} \ \lambda.
\end{equation}

A consequence of \eqref{warunek na miare radialna b}, \eqref{spherical part non-degenerated} and \eqref{nonzerol} is that if $(G,Z)$ is a generating pair, then
\begin{gather}\label{G0 znowu zero}
G(0)=0.
\end{gather}
Indeed, by Proposition \ref{prop wstepny}(A(b)), the jump part of $Z^{G(0)}$ is of finite variation.  Therefore
\begin{align*}
\int_{0}^{+\infty} v \ \nu_{G(0)}(\dd v)&=\int_{\mathbb{R}^d} \langle G(0),y\rangle \nu(\dd y)=\int_{\S^{d-1}}\lambda(d \xi)\int_{0}^{+\infty}\langle G(0),r\xi\rangle\gamma(\dd r)\\
&=\int_{\S^{d-1}}\langle G(0),\xi \rangle \lambda(d \xi)\int_{0}^{+\infty}r \ \gamma(\dd r)<+\infty,
\end{align*}
which, in view of \eqref{warunek na miare radialna b}, \eqref{spherical part non-degenerated} and \eqref{nonzerol} implies \eqref{G0 znowu zero}.

\subsection{Main results}

In this section we prove the following theorem.

\begin{tw} \label{tw_sferyczne}
Let $Z$ be a L\'evy martingale with characteristic triplet $(a,Q,\nu(\dd y))$ such that  $\nu(\dd y)$ 
admits the decomposition \eqref{spherical} with spherical measure $\lambda(\dd \xi)$ satisfying \eqref{spherical part non-degenerated}. Let us also assume that $\gamma(\dd r)$ satisfies \eqref{warunek na miare radialna b} or \eqref{G0 znowu zero} holds. Moreover, let $G:[0,+\infty)\longrightarrow \mathbb{R}^d$ be a continuous function such that  
\begin{gather}\label{istnienie G0}
G_0:=\lim_{x \ra 0+} \frac{G(x)}{|G(x)|},
\end{gather} 
exists. 

Then \eqref{rownanie 2} generates an affine model if and only if $F(x)=ax+b$, $a\in\mathbb{R}, b\geq 0$ and
the measure $\mu(\dd v)$ in Proposition \ref{prop wstepny}(A(c)) is $\alpha$-stable with $\alpha \in (1,2)$.
\end{tw}

The proof of Theorem \ref{tw_sferyczne} is presented in Subsection \ref{proof_of_tw_sf} and is preceded by 
some auxiliary results presented in Subsection \ref{Auxilliary results}.

\vskip 1ex

From Theorem \ref{tw_sferyczne} the following corollary follows.

\begin{cor} \label{cor_sferyczne} Let the assumptions of Theorem \ref{tw_sferyczne} be satisfied. If $(Z,G)$ is a generating pair, then the continuous (Wiener) part of the process $Z^{G(x)}$ vanishes for all $x >0$.
\end{cor}
{\bf Proof:}  It follows from Proposition \ref{prop wstepny}(c) that the Wiener part of $Z^{G(x)}$ satisfies
\begin{equation} \label{warrr4}
\frac{1}{2}\langle QG(x),G(x)\rangle=cx, \ x\geq 0 \ \text{for some} \  c\geq 0.
\end{equation}
Either directly by assumption \eqref{G0 znowu zero} or by assumption \eqref{warunek na miare radialna b} we get that $G(0)=0$. Therefore, by \eqref{war na exp Laplaca} and Theorem \ref{tw_sferyczne}, the Laplace transform of the jump part of $Z$ satisfies 
\begin{equation} \label{warrr3}
J_{\nu}(b G(x))=xJ_{\mu}(b) = \gamma x b^{\alpha}, \quad x\geq 0\ \text{for some}\ \gamma > 0, \alpha \in (1,2).
\end{equation}
Proposition \ref{prop wstepny}(A(a)) guarantees that $\langle G_0,y \rangle \ge 0$ for any $ y \in \text{supp } \nu$ and condition \eqref{spherical part non-degenerated} guarantees that $y \mapsto \langle G_0, y\rangle, y\in  \text{supp} \ \nu$, does not vanish, hence $J_\nu(G_0) > 0$. Consequently, from \eqref{warrr3} we obtain
\[
\lim_{x \ra 0+} \frac{\gamma x}{|G(x)|^{\alpha}} =  \lim_{x \ra  0+} J_{\nu}\rbr{\frac{G(x)}{|G(x)|}} = J_{\nu}\rbr{G_{0}} \in (0, +\ns).
\]
From this, $\lim_{x \ra 0+} {|G(x)|} = 0$ and \eqref{warrr4}  we further have 
\[
\langle Q G_{0},G_0\rangle = \lim_{x \ra 0+} \frac{\langle Q G(x),G(x)\rangle}{|G(x)|^2} =  \lim_{x \ra 0+} \frac{\gamma x}{|G(x)|^{\alpha}}  \frac{2c / \gamma}{|G(x)|^{2 - \alpha}} =
\begin{cases} 
0  \ \ \ \ \ \text{if } c = 0;\\
+ \ns   \ \text{if } c>0.
\end{cases} 
\]
Since  $\langle Q G_{0},G_0\rangle \neq +\ns$, we necessarily have  $c = 0$ which, in view of \eqref{warrr4}, means that the continuous (Wiener) part of $Z^{G(x)}$ vanishes.
\hfill $\square$

\begin{rem} 
A generating pair $(G,Z)$ satisfying assumptions of Theorem \ref{tw_sferyczne} has projections $Z^{G(x)}, x\geq 0$ with the same law as the projections $\tilde{Z}^{\tilde{G}(x)},x\geq 0$, where $\tilde{G}(x)=(C x)^{\frac{1}{\alpha}}, C>0$, and $\tilde{Z}$ is a one-dimensional  $\alpha$-stable process with positive jumps. In view of Proposition \ref{prop wstepny}, the short rates given by \eqref{rownanie 2}
and \eqref{CIR equation geenralized} have the same generator.
\end{rem}

\begin{rem}
In the formulation of Theorem \ref{tw_sferyczne} the assumption \eqref{istnienie G0} can be replaced by the existence 
of the limit $\lim_{x \ra +\ns} \frac{G(x)}{|G(x)|}$. Under the latter condition we were, however, unable to prove Corollary \ref{cor_sferyczne}.
\end{rem}

\begin{rem}
Our proof of Theorem \ref{tw_sferyczne} seems to work for more general measures, namely measures satisfying 
\eqref{spherical} with $\S^{d-1}$ replaced by the boundary $\partial D$ of some convex set $D$ in $\R^d$, containing $0$ in its interior, and the measure $\lambda$ replaced by an appropriate finite measure on $\partial D$.
\end{rem}

\subsection{Auxilliary results}\label{Auxilliary results}

Our first aim is to estimate, for a generating pair $(G,Z)$, the function $J_\nu(b G(x)), b,x\geq 0$ with the use of the function $J_\nu(b G_{0}), b\geq 0$ for $x$ such that $G(x)/\mid G(x)\mid$ is close to $G_0$. 
The solution of this problem is presented in Lemma \ref{lem oszacownaia na funkcje podcalkowa}, Proposition \ref{positivity} and Proposition \ref{sferycznosc}.

Let $\rho(\dd v)$ be a L\'evy measure on $(0,+\infty)$ satisfying
\begin{gather}\label{Ro}
\int_{0}^{+\infty}(v^2\wedge v)\rho(\dd v)<+\infty,
\end{gather}
and 
\begin{gather}\label{funkcja Jot}
J_\rho(z):=\int_{(0,+\infty)}(e^{-zv}-1+zv)\rho(\dd v), \quad z\geq 0,
\end{gather}
Recall, Proposition \ref{bounds_alpha} provides a growth estimation of the function $J_{\rho}$.
The second aim of this section is to provide sufficient conditions for $J_\rho$
to be a power function. This problem is solved in Lemma \ref{lem o charakteryzacji f potegowej} and Lemma \ref{lem Weyl}.

\begin{lem}\label{lem oszacownaia na funkcje podcalkowa}
The function $H:[0,+\infty)\longrightarrow \mathbb{R}$ given by
$$
H(z)=e^{-z}-1+z,
$$
is convex, strictly increasing and 
\begin{gather}\label{min max}
\min\{1,t^2\}\cdot H(z)\leq H(t z)\leq\max\{1,t^2\}\cdot H(z),\quad z\geq 0,t>0.
\end{gather}
\end{lem}
{\bf Proof:} Since $H^\prime(z)=1-e^{-z}$ the monotonicity and convexity of $H$ follows. For $t\geq 1$ it follows from the monotonicity of $H$ that
$$
H(tz)\geq H(z)=\min\{1,t^2\}H(z).
$$
From \eqref{oszH} we obtain
$$
\frac{d}{ds}\ln H(s)=\frac{H^\prime(s)}{H(s)}=\frac{1-e^{-s}}{e^{-s}-1+s}\leq \frac{2}{s}, \quad s>0,
$$
and, consequently, we obtain that for $t\geq 1$: 
$$
\ln H(tz)-\ln H(z)\leq \int_{z}^{tz}\frac{2}{s}ds=\ln t^2.
$$
Thus 
\begin{gather}\label{ktktk}
\min\{1,t^2\}H(z) = H(z)  \le H(tz)\leq t^2 H(z)=\max\{1,t^2\}H(z).
\end{gather}
Using the monotonicity of $H$ and  \eqref{ktktk} we see that for $t\in(0,1)$:
$$
H(tz)\leq H(z)=H\rbr{\frac{1}{t}tz}\leq \frac{1}{t^2}H(tz),
$$
so also for $t\in(0,1)$
$$
\min\{1,t^2\}H(z)=t^2 H(z)\leq H(tz)\leq H(z)=\max\{1,t^2\} H(z).
$$
\hfill $\square$

\begin{cor} \label{lem oszacownaia na wzrost transformaty} Let $\rho(\dd v)$ be a L\'evy measure on $(0,+\infty)$ satisfying $\int_{0}^{+\infty}(v^2\wedge v)\rho(\dd v)<+\infty$. 
It follows from \eqref{min max} and the formula
$$
J_\rho(z):=\int_{(0,+\infty)}H(zv)\rho(\dd v)<+\infty
$$
that the function $J_\rho$ satisfies
\begin{gather}\label{min maxxx}
\min\cbr{1,t^2}\cdot J_\rho(z)\leq J_\rho(t z)\leq\max\cbr{1,t^2}\cdot J_\rho(z),\quad z\geq 0,t>0.
\end{gather}
\end{cor}

\begin{prop} \label{positivity} Let $Z$ be a L\'evy process with characteristic triplet $(a,Q,\nu(\dd y))$.
If \eqref{rownanie 2} generates an affine model and $G_{\ns}$ is an arbitraty limit point of the set
\[
\cbr{\frac{G(x)}{\left|G(x)\right|}:x>0}
\]
then
\[
\nu\cbr{y\in\R^{d}:\left\langle G_{\ns},y\right\rangle <0}=0.
\]
\end{prop}
{\bf Proof:}  Assume that
\[
\nu\cbr{y\in\R^{d}:\left\langle G_{\ns},y\right\rangle <0}=\nu\cbr{y\in\R^{d}\setminus \{0\}:\left\langle G_{\ns},\frac{y}{\left|y\right|}\right\rangle <0}>0.
\]
Then there exists a natural $n$ such that for
\[
V_{n}:=\cbr{y\in\R^{d}\setminus\{0\}:\left\langle G_{\ns},\frac{y}{\left|y\right|}\right\rangle <-\frac{1}{n}}
\]
 one has $\nu\rbr{V_{n}}>0.$

Let $x$ be such that 
\[
\left|\frac{G(x)}{\left|G(x)\right|}-G_{\ns}\right|\le\frac{1}{2n}.
\]
It follows from the Schwarz inequality that for any $y\in\R^{d}$,
\begin{equation}
\left|\left\langle \frac{G(x)}{\left|G(x)\right|},y\right\rangle -\left\langle G_{\ns},y\right\rangle \right|\le\left|\frac{G(x)}{\left|G(x)\right|}-G_{\ns}\right|\left|y\right|\le\frac{1}{2n}\left|y\right|.\label{eq:schwarz}
\end{equation}
Let $y\in V_{n}$. From (\ref{eq:schwarz}) and the definition of
$V_{n}$ we estimate
\[
\left\langle \frac{G(x)}{\left|G(x)\right|},y\right\rangle \le\left\langle G_{\ns},y\right\rangle +\frac{1}{2n}\left|y\right|<-\frac{1}{n}\left|y\right|+\frac{1}{2n}\left|y\right|=-\frac{1}{2n}\left|y\right|<0.
\]
Hence 
\[
\nu\cbr{y\in\R^{d}:\left\langle \frac{G(x)}{\left|G(x)\right|},y\right\rangle <0}\ge\nu\rbr{V_{n}}>0
\]
which is a contradiction with Proposition \ref{prop wstepny}(A(a)).
\hfill$\square$

\begin{prop} \label{sferycznosc} 
Let $(G,Z)$ be a generating pair where $Z$ is a spherical L\'evy process. Assume that $\nu(\dd y)$ has the form \eqref{spherical} and  \eqref{spherical part non-degenerated} holds. Let $G_{\ns}$ be any limit point of the set 
\[
\cbr{\frac{G(x)}{\left|G(x)\right|}:x>0}.
\]
Define
\[
M_{G_{\ns}}(b):=J_{\nu}\rbr{b\cdot G_{\ns}} :=\int_{\S^{d-1}}\int_{0}^{+\ns}H\rbr{b\left\langle G_{\ns},r\cdot \xi\right\rangle }\gamma\rbr{\dd r}\lambda(\dd \xi),
\]
where $H(z):= e^{-z}-1+z$.
There exists a function $\delta: (0,1) \ra (0,+\ns)$ such that for any $\varepsilon_{0}>0$, any $b\ge0$ and $x>0$ such that $\left|\frac{G(x)}{\left|G(x)\right|}-G_{\ns}\right|\le \delta\rbr{\varepsilon_0}$ we have
 
\begin{equation}
\left(1-\varepsilon_{0}\right)M_{G_{\ns}}(b\left|G(x)\right|)\le J_{\nu}\rbr{bG(x)}\le\left(1+\varepsilon_{0}\right)M_{G_{\ns}}(b\left|G(x)\right|).\label{eq:teeza}
\end{equation}
\end{prop}
{\bf Proof:}  Let $\varepsilon\in(0,1)$ be such that 
\begin{equation}
\left(1+\varepsilon\right)^{2}\left(1+\frac{4\varepsilon}{\left(1-\varepsilon\right)^{3}}\right)\le1+\varepsilon_{0},\quad\frac{\left(1-\varepsilon\right)^{2}}{\left(1+\frac{\varepsilon}{1-\varepsilon}\right)}\ge1-\varepsilon_{0}.\label{eq:mniam_mniam_mniam}
\end{equation}
Let us assume that 
\begin{equation} \label{normalisation}
\lambda\cbr{\xi\in\S^{d-1}:\left\langle G_{\ns},\xi\right\rangle >0} = \lambda\rbr{\S^{d-1}}-\lambda\cbr{\xi\in\S^{d-1}:\left\langle G_{\ns},\xi\right\rangle =0}=1,
\end{equation}
(we can assume this, multiplying $\lambda(\dd \xi)$ by a positive constant, provided \[\lambda\cbr{\xi\in\S^{d-1}:\left\langle G_{\ns},\xi\right\rangle >0}>0,\]
otherwise it follows from Proposition \ref{positivity} that we get a
degenerated case \[\lambda\rbr{\S^{d-1}}=\lambda\cbr{\xi\in\S^{d-1}:\left\langle G_{\ns},\xi\right\rangle =0}\] 
where \eqref{spherical part non-degenerated} is broken).
Let $\eta\in(0,1)$ be such that
\begin{equation} \label{normalisation1}
\lambda\cbr{\xi\in\S^{d-1}:0<\left\langle G_{\ns},\xi\right\rangle <\eta}\le\varepsilon.
\end{equation}
Moreover, by Proposition \ref{positivity}, 
\[
\nu\cbr{y\in\R^{d}:\left\langle G_{\ns},y\right\rangle <0}=\lambda\cbr{\xi\in\S^{d-1}:\left\langle G_{\ns},\xi\right\rangle <0}\cdot\gamma\rbr{\R_{+}}=0.
\]
Let us define 
\[
\mathbb{V}_{\eta}=\cbr{\xi\in\S^{d-1}:0<\left\langle G_{\ns},\xi\right\rangle <\eta}.
\]
Let $x$ be such that 
\[
\left|\frac{G(x)}{\left|G(x)\right|}-G_{\ns}\right|\le \delta\rbr{\varepsilon_0} :=  \eta\cdot\varepsilon.
\]
From Lemma \ref{lem oszacownaia na funkcje podcalkowa}, for $b,r\ge0$ and $\xi\in\S^{d-1}$ such that $\left\langle G_{\ns},\xi\right\rangle \in[0,\eta)$
we estimate 
\begin{align*}
  H&\rbr{b\cdot r\left\langle G(x),\xi\right\rangle } \le H\rbr{b\cdot r\left|G(x)\right|\left(\left\langle G_{\ns},\xi\right\rangle +\left|\left\langle \frac{G(x)}{\left|G(x)\right|}-G_{\ns},\xi\right\rangle \right|\right)}\\
 & \le\max\left\{H\rbr{b\cdot r\left|G(x)\right|2\left\langle G_{\ns},\xi\right\rangle },H\rbr{b\cdot r\cdot\left|G(x)\right|2\left|\left\langle \frac{G(x)}{\left|G(x)\right|}-G_{\ns},\xi\right\rangle \right|}\right\}\\
 & \le\max\left\{H\rbr{b\cdot r\left|G(x)\right|2\eta},H\rbr{b\cdot r\cdot\left|G(x)\right|2\eta\cdot\varepsilon}\right\}\\
 & =H\rbr{b\cdot r\left|G(x)\right|2\eta}\\
 & \le4H\rbr{b\cdot r\left|G(x)\right|\eta}.
\end{align*}
Hence 
\begin{align}
\int_{\mathbb{V}_{\eta}}\int_{0}^{+\ns}H\rbr{b\cdot r\left\langle G(x),\xi\right\rangle }\gamma\rbr{\dd r}\lambda(\dd \xi)
 & \le4\int_{\mathbb{V}_{\eta}}\int_{0}^{+\ns}H\rbr{b\cdot r\left|G(x)\right|\eta}\gamma\rbr{\dd r}\lambda(\dd \xi)\nonumber \\
 & \le4\varepsilon\int_{0}^{+\ns}H\rbr{b\cdot r\left|G(x)\right|\eta}\gamma\rbr{\dd r}.\label{eq:osz_kur}
\end{align}
From Lemma \ref{lem oszacownaia na funkcje podcalkowa}, for $b,r\ge0$ and $\xi\in\S^{d-1}$ such that $\left\langle G_{\ns},\xi\right\rangle \in[\eta,1]$, 
we also estimate 
\begin{align}
 H\rbr{b\cdot r\left\langle G(x),\xi\right\rangle } & \le H\rbr{b\cdot r\left|G(x)\right|\left(\left\langle G_{\ns},\xi\right\rangle +\left|\left\langle \frac{G(x)}{\left|G(x)\right|}-G_{\ns},\xi\right\rangle \right|\right)}\nonumber \\
 & \le H\rbr{b\cdot r\left|G(x)\right|\left(\left\langle G_{\ns},\xi\right\rangle +\left\langle G_{\ns},\xi\right\rangle \varepsilon\right)}\nonumber \\
 & \le\left(1+\varepsilon\right)^{2}H\rbr{b\cdot r\left|G(x)\right|\left\langle G_{\ns},\xi\right\rangle },\label{eq:osz_kur_3}
\end{align}
and 
\begin{align}
H\rbr{b\cdot r\left\langle G(x),\xi\right\rangle }\nonumber 
 & \ge H\rbr{b\cdot r\left|G(x)\right|\left(\left\langle G_{\ns},\xi\right\rangle -\left|\left\langle \frac{G(x)}{\left|G(x)\right|}-G_{\ns},\xi\right\rangle \right|\right)}\nonumber \\
 & \ge H\rbr{b\cdot r\left|G(x)\right|\left(\left\langle G_{\ns},\xi\right\rangle -\left\langle G_{\ns},\xi\right\rangle \varepsilon\right)}\nonumber \\
 & \ge\left(1-\varepsilon\right)^{2}H\rbr{b\cdot r\left|G(x)\right|\left\langle G_{\ns},\xi\right\rangle }.\label{eq:osz_kur_2}
\end{align}
Notice that by \eqref{normalisation} and \eqref{normalisation1}, $\lambda\rbr{\mathbb{S}^{d-1} \setminus \mathbb{V}_{\eta}} \ge 1-\varepsilon$. From \eqref{eq:osz_kur_2} and then from  $\lambda\rbr{\mathbb{S}^{d-1}\setminus \mathbb{V}_{\eta}} \ge 1-\varepsilon$ and 
 \eqref{eq:osz_kur} we obtain 
\begin{align}\nonumber
&\int_{\S^{d-1}\setminus\mathbb{V}_{\eta}}\int_{0}^{+\ns}H\rbr{b\cdot r\left\langle G(x),\xi\right\rangle }\gamma\rbr{\dd r}\lambda(\dd \xi)\\ \nonumber
 & \ge\int_{\S^{d-1}\setminus\mathbb{V}_{\eta}}\int_{0}^{+\ns}\left(1-\varepsilon\right)^{2}H\rbr{b\cdot r\left|G(x)\right|\left\langle G_{\ns},\xi\right\rangle }\gamma\rbr{\dd r}\lambda(\dd \xi)\nonumber \\
 & \ge\left(1-\varepsilon\right)^{2}\int_{\S^{d-1}\setminus\mathbb{V}_{\eta}}\int_{0}^{+\ns}H\rbr{b\cdot r\left|G(x)\right|\eta}\gamma\rbr{\dd r}\lambda(\dd \xi)\nonumber \\
 & \ge\left(1-\varepsilon\right)^{2}\left(1-\varepsilon\right)\int_{0}^{+\ns}H\rbr{b\cdot r\left|G(x)\right|\eta}\gamma\rbr{\dd r}\nonumber \\
 & \ge\frac{\left(1-\varepsilon\right)^{3}}{4\varepsilon}\int_{\mathbb{V}_{\eta}}\int_{0}^{+\ns}H\rbr{b\cdot r\left\langle G(x),\xi\right\rangle }\gamma\rbr{\dd r}\lambda(\dd \xi).\label{eq:osz_kur_4}
\end{align}
From (\ref{eq:osz_kur_4}) and (\ref{eq:osz_kur_3}) we estimate
\begin{align*}
J_{\nu}\rbr{bG(x)}= & \int_{\S^{d-1}\setminus\mathbb{V}_{\eta}}\int_{0}^{+\ns}H\rbr{b\cdot r\left\langle G(x),\xi\right\rangle }\gamma\rbr{\dd r}\lambda(\dd \xi)+\int_{\mathbb{V}_{\eta}}\int_{0}^{+\ns}H\rbr{b\cdot r\left\langle G(x),\xi\right\rangle }\gamma\rbr{\dd r}\lambda(\dd \xi)\\
\le & \int_{\S^{d-1}\setminus\mathbb{V}_{\eta}}\int_{0}^{+\ns}H\rbr{b\cdot r\left\langle G(x),\xi\right\rangle }\gamma\rbr{\dd r}\lambda(\dd \xi)\\
 & +\frac{4\varepsilon}{\left(1-\varepsilon\right)^{3}}\int_{\S^{d-1}\setminus\mathbb{V}_{\eta}}\int_{0}^{+\ns}H\rbr{b\cdot r\left\langle G(x),\xi\right\rangle }\gamma\rbr{\dd r}\lambda(\dd \xi)\\
\le & \left(1+\varepsilon\right)^{2}\left(1+\frac{4\varepsilon}{\left(1-\varepsilon\right)^{3}}\right)\int_{\S^{d-1}}\int_{0}^{+\ns}H\rbr{b\cdot r\left|G(x)\right|\left\langle G_{\ns},\xi\right\rangle }\gamma\rbr{\dd r}\lambda(\dd \xi)\\
= & \left(1+\varepsilon\right)^{2}\left(1+\frac{4\varepsilon}{\left(1-\varepsilon\right)^{3}}\right)M_{G_{\ns}}\left(b\cdot r\left|G(x)\right|\right).
\end{align*}
Hence 
\begin{equation}
J_{\nu}\rbr{bG(x)}\le\left(1+\varepsilon\right)^{2}\left(1+\frac{4\varepsilon}{\left(1-\varepsilon\right)^{3}}\right)M_{G_{\ns}}\left(b\cdot r\left|G(x)\right|\right).\label{eq:mniam}
\end{equation}
In order to get the lower bound let us notice that 
\begin{align*}
\int_{\S^{d-1}\setminus\mathbb{V}_{\eta}}\int_{0}^{+\ns}H\rbr{b\cdot r\left|G(x)\right|\left\langle G_{\ns},\xi\right\rangle}\gamma\rbr{\dd r}\lambda(\dd \xi)
 & \ge\int_{\S^{d-1}\setminus\mathbb{V}_{\eta}}\int_{0}^{+\ns}H\rbr{b\cdot r\left|G(x)\right|\eta}\gamma\rbr{\dd r}\lambda(\dd \xi)\\
 & \ge\left(1-\varepsilon\right)\int_{0}^{+\ns}H\rbr{b\cdot r\left|G(x)\right|\eta}\gamma\rbr{\dd r}.
\end{align*}
and 
\begin{align*}
 \int_{\mathbb{V}_{\eta}}\int_{0}^{+\ns}H\rbr{b\cdot r\left|G(x)\right|\left\langle G_{\ns},\xi\right\rangle }\gamma\rbr{\dd r}\lambda(\dd \xi) & \le\int_{\mathbb{V}_{\eta}}\int_{0}^{+\ns}H\rbr{b\cdot r\left|G(x)\right|\eta}\gamma\rbr{\dd r}\lambda(\dd \xi)\\
 & \le\varepsilon\int_{0}^{+\ns}H\rbr{b\cdot r\left|G(x)\right|\eta}\gamma\rbr{\dd r}.
\end{align*}
Hence
\begin{align}
 & \int_{\S^{d-1}\setminus\mathbb{V}_{\eta}}\int_{0}^{+\ns}H\rbr{b\cdot r\left|G(x)\right|\left\langle G_{\ns},\xi\right\rangle }\gamma\rbr{\dd r}\lambda(\dd \xi)\nonumber \\
 & \ge\frac{1-\varepsilon}{\varepsilon}\int_{\mathbb{V}_{\eta}}\int_{0}^{+\ns}H\rbr{b\cdot r\left|G(x)\right|\left\langle G_{\ns},\xi\right\rangle }\gamma\rbr{\dd r}\lambda(\dd \xi),\label{eq:osz_kur_5}
\end{align}
and from this we obtain 
\begin{align}
\int_{\S^{d-1}}\int_{0}^{+\ns}&H\rbr{b\cdot r\left|G(x)\right|\left\langle G_{\ns},\xi\right\rangle }\gamma\rbr{\dd r}\lambda(\dd \xi) =\int_{\S^{d-1}\setminus\mathbb{V}_{\eta}}\int_{0}^{+\ns}H\rbr{b\cdot r\left|G(x)\right|\left\langle G_{\ns},\xi\right\rangle }\gamma\rbr{\dd r}\lambda(\dd \xi) \nonumber \\
 & \quad+\int_{\mathbb{V}_{\eta}}\int_{0}^{+\ns}H\rbr{b\cdot r\left|G(x)\right|\left\langle G_{\ns},\xi\right\rangle }\gamma\rbr{\dd r}\lambda(\dd \xi)\label{eq:osz_kur_5-1} \nonumber \\
 & \le\left(1+\frac{\varepsilon}{1-\varepsilon}\right)\int_{\S^{d-1}\setminus\mathbb{V}_{\eta}}\int_{0}^{+\ns}H\rbr{b\cdot r\left|G(x)\right|\left\langle G_{\ns},\xi\right\rangle }\gamma\rbr{\dd r}\lambda(\dd \xi).
\end{align}
From (\ref{eq:osz_kur_2}) and (\ref{eq:osz_kur_5-1}) we
get 
\begin{align*}
J_{\nu}\rbr{bG(x)}\ge & \int_{\S^{d-1}\setminus\mathbb{V}_{\eta}}\int_{0}^{+\ns}H\rbr{b\cdot r\left\langle G(x),\xi\right\rangle }\gamma\rbr{\dd r}\lambda(\dd \xi)\\
\ge & \left(1-\varepsilon\right)^{2}\int_{\S^{d-1}\setminus\mathbb{V}_{\eta}}\int_{0}^{+\ns}H\rbr{b\cdot r\left|G(x)\right|\left\langle G_{\ns},\xi\right\rangle }\gamma\rbr{\dd r}\lambda(\dd \xi)\\
\ge & \frac{\left(1-\varepsilon\right)^{2}}{\left(1+\frac{\varepsilon}{1-\varepsilon}\right)}\int_{\S^{d-1}}\int_{0}^{+\ns}H\rbr{b\cdot r\left|G(x)\right|\left\langle G_{\ns},\xi\right\rangle }\gamma\rbr{\dd r}\lambda(\dd \xi)\\
= & \frac{\left(1-\varepsilon\right)^{2}}{\left(1+\frac{\varepsilon}{1-\varepsilon}\right)}M_{G_{\ns}}\left(b\cdot r\left|G(x)\right|\right).
\end{align*}
Hence 
\begin{equation}
J_{\nu}\rbr{bG(x)}\ge\frac{\left(1-\varepsilon\right)^{2}}{\left(1+\frac{\varepsilon}{1-\varepsilon}\right)}M_{G_{\ns}}\left(b\cdot r\left|G(x)\right|\right).\label{eq:mniam_mniam}
\end{equation}
Now (\ref{eq:teeza}) follows from (\ref{eq:mniam}), (\ref{eq:mniam_mniam})
and (\ref{eq:mniam_mniam_mniam}).
\hfill$\square$

\begin{lem}\label{lem o charakteryzacji f potegowej}
Let $J_{\rho}$ be given by \eqref{funkcja Jot} with $\rho(\dd v)$ satisfying \eqref{Ro}. Assume that 
\begin{gather}\label{beta war do stabilnosci}
J_\rho(\beta b)=\eta J_{\rho}(b), \quad b\geq 0,
\end{gather}
and 
\begin{gather}\label{gamma war do stabilnosci}
J_\rho(\gamma b)=\theta J_{\rho}(b), \quad b\geq 0,
\end{gather}
for some $\beta>1$, $\gamma>1$ such that $\ln\beta/\ln\gamma\notin \mathbb{Q}$ and $\eta>1$, $\theta>1$. Then
\begin{gather}\label{wzor do udowodnienia}
J_\rho(b)=C b^{\alpha}, b\geq 0,
\end{gather}
for some $C>0$ and $\alpha\in(1,2)$.
\end{lem}
{\bf Proof :} By iterative application of \eqref{beta war do stabilnosci} and \eqref{gamma war do stabilnosci} we see that for any 
$m,n\in\mathbb{N}$
\begin{gather*}
J_{\rho}(\beta^m\gamma^n b)=\eta^m\theta^n J_\rho(b), \quad b\geq 0,
\end{gather*}
which can be written as
\begin{gather}\label{wyjsciowa rownosc w dowodzie}
J_{\rho}(b e^{m \ln\beta+n \ln\gamma})=e^{m\ln\eta+n \ln\theta} J_\rho(b), \quad b\geq 0.
\end{gather}
In Lemma \ref{lem Weyl} below we prove that the set 
$$
D:=\{m \ln\beta -n\ln\gamma; \ m,n\in\mathbb{Z}\}
$$
is dense in $\mathbb{R}$. So, for any $\delta>0$ there exist $m,n\in\mathbb{Z}, m\neq 0,$ such that
\begin{gather}\label{odleglosc 1}
\mid m \ln\beta -n\ln\gamma\mid<\delta,
\end{gather}
and then, by \eqref{min maxxx} and \eqref{wyjsciowa rownosc w dowodzie}, we obtain that
\begin{gather}\label{odleglosc 2}
e^{-2 \delta}\leq \frac{e^{m\ln \eta}}{e^{n\ln \theta}}=\frac{J_\rho(e^{m\ln \beta})}{J_\rho(e^{n\ln \gamma})}\leq e^{2\delta}.
\end{gather}
It follows from \eqref{odleglosc 1} that
$$
\left\vert\frac{\ln\beta}{\ln\gamma}-\frac{n}{m}\right\vert\leq \frac{\delta}{\mid m\mid\ln\gamma},
$$
and from \eqref{odleglosc 2} that
$$
\left\vert\frac{\ln\eta}{\ln\theta}-\frac{n}{m}\right\vert\leq \frac{2\delta}{\mid m\mid\ln\theta}.
$$
Consequently,
$$
\left\vert\frac{\ln\beta}{\ln\gamma}-\frac{\ln\eta}{\ln\theta}\right\vert\leq\frac{\delta}{\mid m\mid\ln\gamma}+
\frac{2\delta}{\mid m\mid\ln\theta}\leq\frac{\delta}{\ln\gamma}+
\frac{2\delta}{\ln\theta}.
$$
Letting $\delta\longrightarrow 0$ yields
$$
\frac{\ln\beta}{\ln\gamma}=\frac{\ln\eta}{\ln\theta}.
$$
Let us define
$$
\alpha:=\frac{\ln\eta}{\ln \beta}=\frac{\ln\theta}{\ln\gamma}>0,
$$
and put $b=1$ in \eqref{wyjsciowa rownosc w dowodzie}. This gives
$$
J_\rho(e^{m \ln\beta+n \ln\gamma})=J_\rho(1)\left(e^{m \ln\beta+n \ln\gamma}\right)^\alpha,
$$
which means that $J_\rho(b)=J_\rho(1)b^\alpha$ for $b$ from the set $e^D$
which is dense in $[0,+\infty)$. As $J_\rho$ is continuous, \eqref{wzor do udowodnienia} follows.
Finally, by Proposition \ref{bounds_alpha} it follows that $\alpha\in(1,2)$. 
\hfill$\square$

\vskip1ex
The following result is strictly related to Weyl's equidistribution theorem, see \cite{Weyl}.

\begin{lem}\label{lem Weyl}
Let $p,q >0$ be such that $p/q\notin\mathbb{Q}$. Let us define the set
$$
G:=\{mp+nq; \quad m,n,=1,2,...\}.
$$
Then for each $\delta>0$ there exists a number $M(\delta)>0$ such that
$$
\forall x\geq M(\delta) \quad  \exists \ g\in G \quad\text{such that} \ \mid x-g\mid\leq \delta.
$$
Moreover, the set
$$
D:=\{mp+nq; \quad m,n\in\mathbb{Z}\},
$$
is dense in $\mathbb{R}$.
\end{lem}
{\bf Proof:} Since $p/q\notin\mathbb{Q}$, at least one of $p,q$, say $q$, is irrational. For simplicity assume that $p=1$ and consider the sequence 
$$
r(jq), j=1,2,... \quad \text{where} \  r(x):=x \ \text{mod} \ 1,
$$
of fractional parts of the numbers $jq, j=1,2,...$\ . Recall, Weyl's equidistribution theorem states that 
\begin{gather}\label{weyl result}
\lim_{N\longrightarrow +\infty}\frac{\sharp\{j\leq N: r(jq)\in[a,b]\}}{N}=b-a
\end{gather}
for any $[a,b]\subseteq [0,1)$ if and only if $q$ is irrational.

For fixed $\delta>0$ and $n$ such that $1/n<\delta$ let us consider a partition of $[0,1)$ of the form
$$
[0,1)=\bigcup_{k=0}^{n-1}A_k, \quad A_k:=[k/n, (k+1)/n).
$$
For a natural number $N$ let us consider the set $R_N:=\{r(jq) : j=1,2,...,N\}$. By \eqref{weyl result}, for each $k=0,1,...,n-1$, there exists $N_k$ such that 
$$
R_{N_k}\cap A_k\neq \emptyset.
$$
Then for $\bar{N}:=\max\{N_0,N_1,...,N_{n-1}\}$ we have
$$
R_{\bar{N}}\cap A_k\neq \emptyset, \quad k=0,1,...,n-1.
$$
Let $M=M(\delta):=\bar{N}q$. Then, for $x\geq M$, there exists a number $N_x\leq \bar{N}$ such that
\begin{gather}\label{oszcownaie reszt}
\mid r(N_x q)-r(x)\mid\leq \frac{1}{n}.
\end{gather}
Then for the number 
$$
g:=\lfloor x \rfloor-\lfloor N_x q\rfloor+N_xg\in G
$$
the following holds
\begin{align*}
\mid x-g\mid&=\mid x-(\lfloor x \rfloor-\lfloor N_x q\rfloor+N_xq)\mid\\[1ex]
&=\mid\lfloor x \rfloor+r(x)-\lfloor x \rfloor+\lfloor N_x q\rfloor-N_xq\mid\\[1ex]
&=\mid r(x)-r(N_xq)\mid\leq 1/n<\delta,
\end{align*}
where the last inequality follows from \eqref{oszcownaie reszt}.

The density of $D$ is an immediate consequence of the first part of the Lemma. Indeed, for $x<M(\delta)$ and $g\in G$ such that $x+g>M(\delta)$ there exists $\tilde{g}\in G$ such that $\mid x+g-\tilde{g}\mid<\delta$.

The general case with $p\neq 1$ can be proven in the same way but requires a generalized version of Weyl's theorem, which says that the numbers $r_p(nq), n=1,2,...$, where $r_p(x):=x \ mod \ p$, are equidistributed on $[0,p)$ if and only if $p/q\notin \mathbb{Q}$. This can be proven by a straightforward modification of the original arguments of Weyl.
\hfill $\square$

\subsection{Proof of Theorem \ref{tw_sferyczne}}  \label{proof_of_tw_sf}
By \eqref{war na exp Laplaca} the Laplace transform $J_\nu$ of the the jump part of $Z$ satisfies
\begin{equation} \label{warrr1}
J_{\nu}(b G(x))=J(b G(0))+xJ_{\mu}(b), \quad b,x\geq 0,
\end{equation}
with ${J}_{\mu}(b)$ given by \eqref{def Jmu i Jnu0}, where $\mu(\dd v)$ is a measure satisfying conditions of Proposition \ref{prop wstepny}(c).
By discussion at the beginning of this section we have $G(0) = 0$, hence  \eqref{warrr1} simpilfes to
\begin{equation} \label{warrr2}
J_{\nu}(b G(x))=xJ_{\mu}(b), \quad x\geq 0,
\end{equation}

From the assumption that supp $\nu$ is not contained in any proper linear subspace of $\R^d$ and \eqref{warrr2}, 
we have that $J_{\nu}(y), J_{\mu}(b)>0$, $G(x) \neq 0$, for $y \in \R^d \setminus \cbr{0}, b>0$, $x>0$.

Let $G_{0} = \lim_{x \ra 0+} \frac{G(x)}{|G(x)|}$. 
From Proposition \ref{sferycznosc} it follows that there exists a function $\delta:(0,+\infty)\rightarrow(0,+\infty)$,
such that for any $\varepsilon>0$ from the inequality 
\[
\left|\frac{G(x)}{\left|G(x)\right|}-G_{0}\right|\le\delta(\varepsilon).
\]
follows that for any $b\ge0$
\[
1-\varepsilon\le\frac{J_{\nu}\rbr{b\frac{G(x)}{\left|G(x)\right|}}}{J_{\nu}\rbr{bG_{0}}}\le1+\varepsilon.
\]
Thus for any $\varepsilon>0$ there exists $m(\varepsilon) >0$,
such that for $x\in \rbr{0, m(\varepsilon)}$
\[
\left|\frac{G(x)}{\left|G(x)\right|}-G_{0}\right|\le\delta(\varepsilon),
\]
and hence for any $b>0$
\[
1-\varepsilon\le\frac{J_{\nu}\rbr{b\frac{G(x)}{\left|G(x)\right|}}}{J_{\nu}\rbr{bG_{0}}}\le1+\varepsilon.
\]
Let us fix $\beta>1$ and take $x_{1},x_{2}$ satisfying
$0 < x_{1}\le x_{2} < m(\varepsilon)$, $\beta\left|G(x_{1})\right|=\left|G(x_{2})\right| > 0$
(from the continuity of $G$  it follows that such $x_{1}$ and $x_{2}$ exist).
Then for any $b>0$ and $i=1,2$, by \eqref{warrr2}, 
\[
1-\varepsilon\le\frac{J_{\nu}\rbr{b\frac{G(x_{i})}{\left|G(x_{i})\right|}}}{J_{\nu}\rbr{bG_{0}}}=\frac{x_{i}J_{\mu}\rbr{\frac{b}{\left|G(x_{i})\right|}}}{J_{\nu}\rbr{bG_{0}}}\le1+\varepsilon.
\]
Hence for any $b>0$, taking $\tilde{b}=\beta\left|G(x_{1})\right|b$ we get
\[
\frac{1-\varepsilon}{1+\varepsilon}\cdot \frac{x_{2}}{x_{1}}\le\frac{J_{\mu}\rbr{\frac{\tilde{b}}{\left|G(x_{1})\right|}}}{J_{\mu}\rbr{\frac{\tilde{b}}{\left|G(x_{2})\right|}}}=\frac{J_{\mu}\rbr{\beta b}}{J_{\mu}\rbr b}\le\frac{1+\varepsilon}{1-\varepsilon}\cdot \frac{x_{2}}{x_{1}}
\]
which yields
\[
\frac{1-\varepsilon}{1+\varepsilon}\cdot \frac{J_{\mu}\rbr{\beta b}}{J_{\mu}\rbr b}\le\frac{x_{2}}{x_{1}}\le\frac{1+\varepsilon}{1-\varepsilon}\cdot \frac{J_{\mu}\rbr{\beta b}}{J_{\mu}\rbr b}.
\]
Since $\varepsilon>0$ is arbitrary, taking $\varepsilon\rightarrow0$
and $x_{1},x_{2}$ satisfying
$0 < x_{1}\le x_{2} <m(\varepsilon)$, $\beta\left|G(x_{1})\right|=\left|G(x_{2})\right|$ we obtain that 
\[
\lim_{\varepsilon\rightarrow0}\frac{x_{2}}{x_{1}}=\eta,
\]
where $\eta = J_{\mu}\rbr{\beta b}/J_{\mu}\rbr b > 1$ is independent from $b>0$.
Hence, for all $b\ge0$ we have 
\[
J_{\mu}\rbr{\beta b}=\eta J_{\mu}\rbr b.
\]

Similarly, take $\gamma>1$ such that $\ln\beta/\ln\gamma\notin \mathbb{Q}$. Reasoning similarly as before we get that there exists $\theta>1$, such that for all $b\ge0$ we have
\[
J_{\mu}\rbr{\gamma b}=\theta J_{\mu}\rbr b.
\]
Now the thesis follows from Lemma \ref{lem o charakteryzacji f potegowej} and the one to one correspondence
between Laplace transforms and measures on $[0,+\ns)$, see \cite{{Feller}} p. 233.
\hfill $\square$

\section{Appendix}

\noindent{\bf Proof of Proposition \ref{prop wstepny}:}
$(A)$ It was shown in \cite [Theorem 5.3]{FilipovicATS} that the generator of a general positive Markovian short rate generating  an affine model is of the form
\begin{align}\label{generator Filipovica}
\mathcal{A}f(x)=&c x f^{\prime\prime}(x)+(\beta x+\gamma)f^\prime(x)\\[1ex] \nonumber
&+\int_{(0,+\infty)}\Big(f(x+y)-f(x)-f^\prime(x)(1\wedge y)\Big)(m(\dd y)+x\mu(\dd y)), \quad x\geq 0,
\end{align}
for $f\in\mathcal{L}(\Lambda)\cup C_c^2(\mathbb{R}_{+})$, where 
$\mathcal{L}(\Lambda)$ is the linear hull of $\Lambda:=\{f_\lambda:=e^{-\lambda x}, \lambda\in(0,+\infty)\}$
and $C_c^2(\mathbb{R}_{+})$ stands for the set of twice continuously differentiable functions with compact support in $[0,+\infty)$. 
Above $c, \gamma\geq 0$, $\beta\in\mathbb{R}$ and $m(\dd y)$, $\mu(\dd y)$ are nonnegative Borel measures on $(0,+\infty)$ satisfying
\begin{gather}\label{warunki na iary Filipovica}
\int_{(0,+\infty)}(1\wedge y)m(\dd y)+\int_{(0,+\infty)}(1\wedge y^2)\mu(\dd y)<+\infty.
\end{gather}

The generator of the short rate process given by \eqref{rownanie 2} equals
\begin{align*}
 \mathcal{A}_{R}f(x) =& f^\prime(x)F(x)+\frac{1}{2}f^{\prime\prime}(x)\langle QG(x),G(x)\rangle \\
& +\int_{\mathbb{R}^d}\Big(f(x+\langle G(x),y\rangle)-f(x)-f^\prime(x)\langle G(x),y\rangle\Big)\nu(\dd y) \\
 = & f^\prime(x)F(x)+\frac{1}{2}f^{\prime\prime}(x)\langle QG(x),G(x)\rangle \\
&+\int_{\mathbb{R}}\Big(f(x+v)-f(x)-f^\prime(x)v\Big)\nu_{G(x)}(\dd v)
\end{align*}
where  $f$ is a bounded, twice continuously differentiable function. 

By Proposition \ref{prop o skokach z dodatniosci rozwiazania} below, the support of the measure $\nu_{G(x)}$ is contained in $[-x,+\ns)$, thus it follows that  
\begin{align}\label{Generatorr R gen}\nonumber
\mathcal{A}_{R}f(x) = &f^\prime(x)F(x)+\frac{1}{2}f^{\prime\prime}(x)\langle QG(x),G(x)\rangle \\ \nonumber
&+\int_{(0, +\ns)}\Big(f(x+v)-f(x)-f^\prime(x)(1\wedge v) \Big)\nu_{G(x)}(\dd v)\\ \nonumber
&+f^\prime(x)\int_{(0, +\ns)}\Big((1\wedge v)-v\Big)\nu_{G(x)}(\dd v)\\[1ex] \nonumber 
&+\int_{(-\ns, 0)}\Big(f(x+v)-f(x)-f^\prime(x)v \Big)\nu_{G(x)}(\dd v)\\ \nonumber
= & \frac{1}{2}f^{\prime\prime}(x)\langle QG(x),G(x)\rangle + f^\prime(x)\sbr{F(x)+\int_{(1,+\ns)}\Big(1- v\Big)\nu_{G(x)}(\dd v)} \\[1ex] \nonumber
&+\int_{(0, +\ns)}\Big(f(x+v)-f(x)-f^\prime(x)(1\wedge v) \Big)\nu_{G(x)}(\dd v)\\ 
&+\int_{[-x, 0)}\Big(f(x+v)-f(x)-f^\prime(x)v \Big)\nu_{G(x)}(\dd v).  
\end{align}

Comparing \eqref{Generatorr R gen} with \eqref{generator Filipovica} applied to a function $f_{\lambda}$ with $\lambda >0$ such that $f_{\lambda}(x) = e^{-\lambda x}$ for $x \ge 0$, we get

\begin{align}
&c x \lambda^2  - (\beta x+\gamma)  \lambda  \nonumber \\[1ex] \nonumber
&+ \int_{(0,+\infty)}\Big(e^{-\lambda y}-1+\lambda (1\wedge y)\Big)(m(\dd y)+x\mu(\dd y)) \\[1ex] \nonumber
& -  \frac{1}{2} \lambda^2 \langle QG(x),G(x)\rangle + \sbr{F(x)+\int_{(1,+\ns)}\Big(1- v\Big)\nu_{G(x)}(\dd v)} \lambda  
\\[1ex] \nonumber
& - \int_{(0,+\infty)}\Big(e^{-\lambda v}-1+\lambda (1\wedge v)\Big)\nu_{G(x)}(\dd v) \\[1ex] \label{uaua}
& =  \int_{[-x, 0)}\Big(e^{-\lambda v}-1+\lambda  v \Big)\nu_{G(x)}(\dd v) , \quad \lambda >0, x \geq 0.
\end{align}
Comparing the left and the right sides of \eqref{uaua} we see that the left side grows no faster than a quadratic polynomial of $\lambda$ while the right side grows faster that $d e^{\lambda y}$ for some $d, y >0$, unless  the support of the measure $\nu_{G(x)}(\dd v)$ is contained in $[0,+\ns)$. It follows that $\nu_{G(x)}(\dd v)$ is concentrated on $[0,+\infty)$, hence $(a)$ follows, and 

\begin{align}
&c x \lambda^2  - (\beta x+\gamma)  \lambda  \nonumber \\[1ex] \nonumber
& -  \frac{1}{2} \lambda^2 \langle QG(x),G(x)\rangle + \sbr{F(x)+\int_{(1,+\ns)}\Big(1- v\Big)\nu_{G(x)}(\dd v)} \lambda  
\\[1ex] \label{uaua1}
& = \int_{(0,+\infty)}\Big(e^{-\lambda y}-1+\lambda (1\wedge y)\Big)\rbr{\nu_{G(x)}(\dd y) - m(\dd y)- x\mu(\dd y)}, \quad \lambda >0, x \geq 0.
\end{align}
Dividing both sides of the last equality by $\lambda^2$ and using the  estimate  
$$\frac{e^{-\lambda y}-1+\lambda (1\wedge y)}{\lambda^2} \le \rbr{\frac{1}{2} y^2} \wedge \rbr{\frac{e^{-\lambda}-1+\lambda}{\lambda^2}}$$ 
we get that that the left side of \eqref{uaua1} converges to $c x -  \frac{1}{2} \langle QG(x),G(x)\rangle$ as $\lambda\rightarrow +\infty$, while the right side converges to  $0$. This yields \eqref{mult. CIR condition}, i.e. 
\begin{align}\label{W1}
c x=& \frac{1}{2}\langle QG(x),G(x)\rangle,\quad x\geq 0.
\end{align}
Next, fixing $x\geq 0$ and comparing \eqref{Generatorr R gen} with \eqref{generator Filipovica} applied to a function from the domains of both generators and such that $f(x) = f'(x) = f ''(x) =0$ we get 
\[
\int_{(0,+\infty)} f(x+y) (m(\dd y)+x\mu(\dd y)) =  \int_{(0, +\ns)} f(x+v) \nu_{G(x)}(\dd v) 
\]
for any such a function,
which yields 
\begin{gather}\label{rozklad nuG na sume}
\nu_{G(x)}(\dd v)\mid_{(0,+\infty)}=m(\dd v)+x\mu(\dd v), \quad x\geq 0.
\end{gather} 
This implies also
\begin{align}
\label{W2}
\beta x+\gamma=&F(x)+\int_{(1,+\ns)}\Big(1 - v\Big)\nu_{G(x)}(\dd v),\quad x\geq 0.
\end{align}

$(b)$ Setting $x=0$ in \eqref{rozklad nuG na sume} yields 
\begin{gather}\label{nu G0}
\nu_{G(0)}(\dd v)\mid_{(0,+\infty)}=m(\dd v).
\end{gather}
To prove \eqref{nu G0 finite variation}, by \eqref{warunki na iary Filipovica} and \eqref{nu G0}, we need to show that
\begin{gather}\label{cocococ}
\int_{(1,+\infty)}v\nu_{G(0)}(\dd v)<+\infty.
\end{gather}
It is true if $G(0)=0$ and for $G(0)\neq 0$ the following estimate holds
\begin{align*}
\int_{(1,+\infty)}v\nu_{G(0)}(\dd v)&=\int_{\mathbb{R}^d}\langle G(0),y\rangle\mathbf{1}_{[1,+\infty)}(\langle G(0),y\rangle)\nu(\dd y)\\[1ex]
&\leq \mid G(0)\mid\int_{\mathbb{R}^d}\mid y\mid\mathbf{1}_{[1/\mid G(0)\mid,+\infty)}(\mid y\mid)\nu(\dd y),
\end{align*}
and \eqref{cocococ} follows from \eqref{warunek na miare Levyego 2}.

$(c)$   \eqref{rozklad nu G(x)} follows from \eqref{rozklad nuG na sume} and \eqref{nu G0}.
To prove \eqref{war calkowe na mu} we use \eqref{rozklad nu G(x)},  \eqref{nu G0 finite variation}
and the following estimate for $x\geq 0$:
\begin{align*}
\int_{0}^{+\infty}(v^2\wedge v)\nu_{G(x)}(\dd v)&=\int_{\mathbb{R}^d}(\mid\langle G(x),y\rangle\mid^2\wedge\langle G(x),y\rangle)\nu(\dd y)\\[1ex]
&\leq \Big(\mid G(x)\mid^2\vee \mid G(x)\mid\Big)\int_{\mathbb{R}^d}(\mid y\mid^2\wedge \mid y\mid)\nu(\dd y)<+\infty,
\end{align*}
In the last line we used \eqref{warunek na miare Levyego 1} and \eqref{warunek na miare Levyego 2}.

$(d)$ It follows from \eqref{W2} and \eqref{rozklad nu G(x)} that 
\begin{align*}
\beta x+\gamma&=F(x)+\int_{(1,+\infty)}(1 - v)\nu_{G(x)}(\dd v)\\[1ex]
&=F(x)+\int_{(1,+\infty)}(1-v)\nu_{G(0)}(\dd v)+x\int_{(1,+\infty)}(1-v)\mu(\dd v), \quad x\geq 0.
\end{align*}
Consequently, \eqref{linear drift}
follows with
$$
a:=\Big(\beta-\int_{(1,+\infty)}(1-v)\mu(\dd v)\Big), \ b:=\Big(\gamma-\int_{(1,+\infty)}(1-v)\nu_{G(0)}(\dd v)\Big),
$$
and $b\geq \int_{(1,+\infty)}(v-1)\nu_{G(0)}(\dd v)$ because $\gamma\geq 0$.

$(B)$ We use \eqref{W2},  \eqref{linear drift} and \eqref{rozklad nuG na sume} to write \eqref{generator Filipovica} in the form
\begin{align*}
\mathcal{A}f(x)=cx f^{\prime\prime}(x)&+\Big[ax +b+\int_{(1,+\infty)}(1 -v)\nu_{G(x)}(\dd v)\Big]f^{\prime}(x)\\[1ex]
&+\int_{(0,+\infty)}[f(x+v)-f(x)-f^{\prime}(x)(1\wedge v)]\nu_{G(x)}(\dd v)\}.
\end{align*}
In view of \eqref{rozklad nuG na sume} and \eqref{nu G0} we see that
\eqref{generator R w tw} is true.

\begin{prop}\label{prop o skokach z dodatniosci rozwiazania}
Let $G:[0,+\infty)\rightarrow \mathbb{R}^d$ be continuous. If the equation \eqref{rownanie 2} has a non-negative strong solution for any initial condition $R(0)=x\geq 0$, then
\begin{gather}\label{ograniczenia skokow}
\forall x\geq 0 \quad \nu{\{y\in\mathbb{R}^d: x+\langle G(x),y\rangle<0\}}=0.
\end{gather}
In particular, the support of the measure $\nu_{G(x)}(\dd v)$ is contained in $[-x,+\infty)$.
\end{prop}
{\bf Proof:} Let us assume to the contrary, that for some $x\geq 0$
$$
\nu{\{y\in\mathbb{R}^d: x+\langle G(x),y\rangle<0\}}>0.
$$
Then there exists $c>0$ such that
$$
\nu{\{y\in\mathbb{R}^d: x+\langle G(x),y\rangle<-c\}}>0.
$$
Let $A\subseteq \{y\in\mathbb{R}^d: x+\langle G(x),y\rangle<-c\}$ be a Borel set separated from zero. By the continuity of $G$ we have that for some $\varepsilon>0$:
\begin{gather}\label{pani w szpileczkach}
\tilde{x}+\langle G(\tilde{x}),y\rangle<-\frac{c}{2}, \quad  \tilde{x}\in[(x-\varepsilon)\vee 0,x+\varepsilon],\quad y\in A.
\end{gather}
Let $Z^2$ be a L\'evy processes with characteristics $(0,0,\nu^2(dy))$, where $\nu^2(dy):=\mathbf{1}_{A}(y)\nu(dy)$ and $Z^1$ be defined by $Z(t)=Z^1(t)+Z^2(t)$. Then $Z^1, Z^2$ are independent and $Z^2$ is a compound Poisson process. Let us consider 
the following equations
\begin{gather*}
dR(t)=F(R(t))dt+\langle G(R(t-)),dZ(t)\rangle, \quad R(0)=x,\\[1ex]
dR^1(t)=F(R^1(t))dt+\langle G(R^1(t-)),dZ^1(t)\rangle, \quad R^1(0)=x.
\end{gather*}
For the exit time $\tau_1$ of $R^1$ from the set $[(x-\varepsilon)\vee 0,x+\varepsilon]$ and the first jump time $\tau_2$ of $Z^2$ we can find
$T>0$ such that $\mathbb{P}(\tau_1>T, \tau_2<T)=\mathbb{P}(\tau_1>T)\mathbb{P}(\tau_2<T)>0$. On the set $\{\tau_1>T, \tau_2<T\}$ we have $R(\tau_2-)=R^1(\tau_2-)$ and therefore 
$$
R(\tau_2)=R^1(\tau_2-)+\langle G(R^1(\tau_2-)),\triangle Z^2(\tau_2)\rangle<-\frac{c}{2}.
$$
In the last inequality we used \eqref{pani w szpileczkach}. This contradicts the positivity of $R$. \hfill $\square$
\vskip2ex
\noindent
{\bf Proof of Proposition \ref{prop2}:} The HJM condition for affine models  takes the form
\begin{gather}\label{HJM dla affine}
J_Z(B(v)G(x))=-A^\prime(v)-[B^\prime(v)-1]x+B(v)F(x), \quad v,x\geq 0,
\end{gather}
for details see Proposition 3.2 in \cite{BarskiZabczykCIR}. Using \eqref{linear drift} we obtain
\begin{gather}\label{slonie dwa}
J_Z(B(v)G(x))=-A^\prime(v)+bB(v)+[aB(v)-B^\prime(v)+1]x, \quad v,x\geq 0.
\end{gather}
Setting $x=0$ and using \eqref{war na exp Laplaca} yields
$$
-A^\prime(v)+bB(v)=J_Z(B(v)G(0))=J_{\nu_{G(0)}}(B(v)),
$$
which is \eqref{A wzor od B}. It follows from \eqref{slonie dwa} that
\begin{gather}\label{wzwzorekk}
J_Z(B(v)G(x))=J_Z(B(v)G(0))+[aB(v)-B^\prime(v)+1]x, \quad v,x\geq 0.
\end{gather}
Using again \eqref{war na exp Laplaca} we obtain
\begin{align*}
\frac{1}{2}B^2(v)cx+J_{\nu_{G(0)}}(B(v))+xJ_\mu(B(v))&=J_{\nu_{G(0)}}(B(v))\\[1ex]
&+[aB(v)-B^\prime(v)+1]x, \quad v,x\geq 0.
\end{align*}
Consequently, 
$$
\frac{1}{2}B^2(v)c+J_\mu(B(v))=aB(v)-B^\prime(v)+1, \quad v,x\geq 0,
$$
which finally yields \eqref{rr na B}.
\hfill$\square$

\end{document}